\documentclass[11pt,reqno]{amsart}
\usepackage{amsmath,amstext,amsfonts,amsthm,amssymb}
\usepackage{eucal}
\usepackage{graphicx}
\renewcommand{\subsubsection}[1]{\addtocounter{subsubsection}{1}
{\ \\[3pt]\bf \thesubsubsection. \  #1} }

\swapnumbers

{  \theoremstyle{definition}

}
\numberwithin{equation}{subsection}
\newcommand{\Spec}{\operatorname{Spec}}

\newcommand{\Higgs}{\operatorname{Higgs}}

\newcommand{\iso}{\overset{\sim}{\longrightarrow}}
\newcommand{\isom}{\overset{\sim}{=}}
\newcommand{\lra}{\longrightarrow}
\newcommand{\ra}{\rightarrow}

\newcommand{\bea}{\begin{eqnarray*}}
\newcommand{\eea}{\end{eqnarray*}}
\newcommand{\bean}{\begin{eqnarray}}
\newcommand{\eean}{\end{eqnarray}}

\newcommand{\hfg}{\hat{\frak{g}}}

\newcommand{\hX}{\hat{X}}

\newcommand{\fg}{\mathfrak g}

\newcommand{\ft}{\mathfrak t}

\newcommand{\CA}{\mathcal{A}}
\newcommand{\CB}{\mathcal{B}}

\newcommand{\CD}{\mathcal{D}}
\newcommand{\CF}{\mathcal{F}}

\newcommand{\CO}{\mathcal{O}}

\newcommand{\CT}{\mathcal{T}}
\newcommand{\CV}{\mathcal{V}}

\newcommand{\CDO}{\mathcal{CDO}}
\newcommand{\TDO}{\mathcal{TDO}}

\newcommand{\BC}{\mathbb{C}}

\newcommand{\BP}{\mathbb{P}}

\newcommand{\BR}{\mathbb{R}}

\newcommand{\nc}{\newcommand}

\nc{\Id}{\text{Id}}
\nc{\la}{\lambda}

\begin{document}


\newcommand{\Real}{\mathbb R}
\newcommand{\HH}{\mathbb H}
\newcommand{\QQ}{\mathbb Q}
\newcommand{\ZZ}{\mathbb Z}
\newcommand{\LL}{\mathbb L}
\newcommand{\VV}{\mathbb V}
\newcommand{\MM}{\mathbb M}
\newcommand{\PP}{\mathbb P}
\newcommand{\RR}{\mathbb R}

\newcommand{\cA}{\mathcal{A}}
\newcommand{\cB}{\mathcal{B}}
\newcommand{\cC}{\mathcal{C}}
\newcommand{\cext}{\cC ext}
\newcommand{\cD}{{\mathcal{D}}}
\newcommand{\cE}{{\mathcal{E}}}

\newcommand{\cF}{{\mathcal{F}}}
\newcommand{\cG}{{\mathcal{G}}}
\newcommand{\cH}{{\mathcal{H}}}
\newcommand{\cJ}{{\mathcal{J}}}

\newcommand{\cL}{{\mathcal{L}}}
\newcommand{\cM}{{\mathcal{M}}}
\newcommand{\cN}{\mathcal{N}}
\newcommand{\cO}{{\mathcal{O}}}
\newcommand{\cP}{{\mathcal{P}}}
\newcommand{\cQ}{{\mathcal{Q}}}
\newcommand{\cR}{{\mathcal{R}}}
\newcommand{\cS}{{\mathcal{S}}}
\newcommand{\cT}{{\mathcal{T}}}
\newcommand{\cV}{{\mathcal{V}}}
\newcommand{\cW}{{\mathcal{W}}}
\newcommand{\vext}{\cV ext}

\newcommand{\cZ}{{\mathcal{Z}}}

\title{chiral differential operators on abelian varieties}

\author{Fyodor Malikov, Vadim Schechtman}
\maketitle



\bigskip\bigskip
\begin{center}
{\em To Igor Frenkel on his 60th birthday}
\end{center}
\bigskip\bigskip

\section{introduction}
This paper consists of two largely independent parts; they have in common an underlying geometric object, a torus, and an underlying algebro-geometric
object, an algebra of chiral differential operators.

\subsection{ }
In the Part I (Sections 2, 3) we describe a way to get from an algebra of chiral differential 
operators (cdo) on an abelian variety a cdo (actually a family of cdo's) on the dual variety. 
 
Let $X$ be an abelian variety (everything in this note will be over $\BC$), and $\hat X$ be its dual variety. Set 
$\fg = Lie(X) = H^0(X,\CT_X),\ \hfg = Lie(\hX) = H^1(X,\CO_X)$. Let $\CA$ be an algebra of twisted 
differential operators (a tdo) over $X$, cf. [BB], \S 2; its isomorphism class is an element of
$H^1(X,\Omega_X^1 \lra \Omega_X^{2,cl})$. Let $c(\CA)$ be its image under the natural map 
$$
H^1(X,\Omega_X^1 \lra \Omega_X^{2,cl}) \lra H^1(X,\Omega_X^1).
$$ 

We can identify
$$
H^1(X,\Omega_X^1) \iso H^0(X,\Omega_X^1) \otimes  H^1(X,\CO_X) \iso Hom(\fg,\hfg)
\eqno{(1.1.1)}
$$
cf. [M]. Let us call a class $c\in H^1(X,\Omega_X^1)$ {\it non-degenerate} 
if $c$ considered as a map $\hfg \lra \fg$ is an isomorphism. We denote by 
$H^1(X,\Omega_X^1)_{nd}\subset H^1(X,\Omega_X^1)$ the (Zarisky open) subspace 
of nondegenerate classes.  

Following [PR] let us call a tdo 
$\CA$ non-degenerate if $c(\CA)$ is nondegenerate. For a non-degenerate $\CA$ Polishchuk and Rothstein define 
a tdo $\Phi(\CA)$ on $\hX$, 
called the {\it Fourier - Mukai transform} of $\CA$, with 
$$
c(\Phi(\CA)) = - c(\CA)^{-1}
\eqno{(1.1.2)}
$$
\subsection{ }
 Let $\CDO(X)$ denote the groupoid of {\it chiral} differential operators (cdo) on $X$, cf. [GMS1] and 2.6 below. 
In \S 2 of this note we construct for each $\mu\in H^1(X,\Omega_X^1)_{nd}$ an equivalence of groupoids
$$
F_\mu:\ \CDO(X) \iso \CDO(\hX), 
\eqno{(1.2.1)}
$$
cf. 2.8.    

Our construction has an elementary linear algebra avatar. Let $V, V'$ be 
two finite dimensional vector spaces of equal dimensions, $M = Isom(V,V')\subset Hom(V,V')$ the variety 
of invertible linear maps $V\iso V'$, $M' = Isom(V',V)$. Define a map 
$$
F:\ M \lra M',\ F(A) = - A^{-1}
\eqno{(1.2.2)}
$$
This is a linear algebra analogue of the Fourier transform (1.1.2). We can identify the total space $TM$ 
of the tangent bundle to $M$ with $M\times Hom(V,V')$. The differential of 
(1.2.1) 
$$
dF:\ TM = M\times Hom(V,V') \lra TM = M'\times Hom(V',V)
$$ 
acts as 
$$
dF(A,B) = ( - A^{-1}, A^{-1}BA^{-1})
\eqno{(1.2.3)}
$$
Our a Fourier - Mukai transform for cdo's (1.2.1) is an analog of (1.2.3). 
The role of $A$ will be played by a non-degenerate class $\mu$ and the role of $B$ by a cdo, cf. (2.8.3). 

Our construction is of cohomological nature and uses the Hodge theory (which in the case 
of abelian varieties is completely algebraic). The same construction gives 
rise to the equivalence of groupoids of usual tdo's on $X$ and on $\hX$; presumably 
the answer coincides with the Polishchuk - Rothstein FM transform; this way 
we get an alternative definition of their transformation.

In sect.~3 we make some remarks on the quasiclassical limit of this construction. 

\subsection{ }  Part II of this paper occupies sections 4--7 and we begin with a disclaimer: there is nothing new in these sections. 
Why write them then? A forgiving reader might consider three points. First,  sections 4--6 contain what they call a rigorous derivation
of the quantum
$\sigma$-model on the torus. Second, we explain exactly how the antiholomorphic part of the theory enters the game -- the subject
that seem to have avoided much of the mathematical literature on the subject. Third, this derivation is  elementary -- as elementary as the derivation of the quantum harmonic  oscillator that can be found in any 
quantum mechanics textbook, e.g. [LL].  In fact, we find it convenient to begin with reminding the reader of the latter model.

Consider the Lagrangian $\cL= 1/2((\partial_\tau x)^2-x^2)d\tau$.  The equation of motion is: $\partial_\tau^2x+x=0$. The variational form
$\gamma=\partial_\tau x\delta x$. It follows that the phase space is $\BR\times\BR$ with coordinates $p=\partial_\tau x$ and $x$, and the Poisson bracket $\{p,x\}=1$.
The vector field $\partial_\tau$ is a symmetry of $\cL$ and the corresponding integral of motion is $H_{\partial _\tau}=1/2(p^2+x^2)$. 

Upon quantization, $C^\infty(\BR\times\BR)$ becomes the algebra of differential operators on the line, $\cD_\BR$, the Hamiltonian $H_{\partial_ \tau}=1/2(-\partial_x^2+x^2)$.
What is the space of states? The operator $H_{\partial_\tau}$ must be diagonalizable with spectrum bounded from below. This does not leave us much room for maneuver.
We notice that up to a constant summand $H_{\partial_\tau}$ is $(1/2)a_-a_+$, with $a_-=x+\partial_x$, $a_+=x-\partial_x$ so that $[H_{\partial_\tau},a_{\pm}]=\pm a_{\pm}$.
We obtain a 3-dimensional Heisenberg algebra $\cH$ spanned by $a_\pm$ and 1 that diagonalizes the Hamiltonian. The space of states is then an induced $\cH$-module spanned by $v$ s.t. $a_-v=0$.
Furthermore, this space has a functional realization as a linear span of $\{L_n(x)e^{-x^2/2}\}$, where $L_n(x)$ is the Laguerre polynomial.

Things are very much the same for the $\sigma$-model on the torus except that now  the phase space is the space of jets of an appropriate trivial bundle over the circle $S^1$.
The space of functions, instead of being a Poisson algebra, is a vertex Poisson or rather a coisson algebra [BD].  The coisson bracket gives a usual Lie bracket on quantities
that can symbolically be written as $\int_{S^1}e^{im\sigma}a(\sigma)d\sigma$, $a(\sigma)$ being an element of the coisson algebra. The 
Hamiltonian $\int_{S^1}H_{\partial_\tau}$
splits as $\int_{S^1}H_{\partial_\tau}=\int H_{\partial_z}+\int H_{\partial_{\bar{z}}}$, where $H_{\partial_z}$, $H_{\partial_{\bar{z}}}$ are integrals of motions obtained by
lifting the indicated vector fields.  Diagonalizing $\int_{S^1}H_{\partial_\tau}$  resembles very much the oscillator case except that
now $\int_{S^1}H_{\partial_z}$ and $\int_{S^1}H_{\partial_{\bar{z}}}$ are diagonalized separately by 2 commuting with each other infinite dimensional Heisenberg algebras; the former  by operators of the form $\int e^{ im\sigma} u(\sigma)d\sigma$, the latter by $\int e^{ im\sigma} \bar{u}(\sigma)d\sigma$.  The role of $a_-$ (i.e., the subalgebra acting
locally nilpotently) in the former case is played by
 $\int e^{im\sigma} u(\sigma)d\sigma$  and in the latter by  $\int e^{- im\sigma} \bar{u}(\sigma)d\sigma$ with $m>0$. To compensate for the sign difference,
 one continues analytically in the former case by setting $z=e^{i\sigma}$ and in the latter by setting $\bar{z}=e^{-i\sigma}$. This leads
 to interpreting $u$ as   holomorphic and  $\bar{u}$  as antiholomorphic. The space of states then becomes a direct sum of tensor products of holomorphic and antihomorphic
factors. It is not a chiral algebra (just as the space of states of a quantum oscillator is not an associative algebra); rather it possesses what Kapustin and Orlov called a vertex algebra structure in [KO].

A similar analysis applies to the case where a torus is replaced with a simple compact Lie group $G$. This case is more interesting: the underlying coisson algebra is one of chiral differential operators twisted by a global 3-form [GMS2], and the 2 algebras diagonalizing the Hamiltonian are the lifts of the action of $\fg=Lie(G)$ on the left and on the right.
 We explain this very briefly in sect.~7.


\subsection{ }
 We  are  grateful to A.Polishchuk,  A.Beilinson, and A.Kapustin for numerous illuminating discussions.
Much of this work was done in the inspiring 
atmosphere of Max-Planck-Institut f\"ur Mathematik in summer  of 2011.  F.M. was partially supported by an NSF grant.

\bigskip\bigskip 

\begin{center} {\bf PART I: CHIRAL FOURIER-MUKAI TRANSFORM}
\end{center}

\section{a fourier-mukai transform for cdo's}


\subsection{Lemma}
 {\it Let $Y$ be a topological space, 
$$
\CF^\cdot:\ 0 \lra \CF^0 \lra \CF^1 \lra \ldots
$$
a complex of sheaves on $Y$. Consider a subcomplex 
$$
\tau_{\leq 1}\CF^\cdot:\ 0 \lra \CF^0 \lra \CF^{1,cl} \lra 0
$$
where $\CF^{1,cl} = \text{Ker} (\CF^1 \lra \CF^2)$. 
 
The induced map in cohomology 
$$
H^i(Y,\tau_{\leq 1}\CF^\cdot) \lra H^i(Y,\CF^\cdot)
$$
is an isomorphism 
for $i = 0, 1$ and a monomorphism for $i=2$.}

(We are obliged to H.Esnault for this elementary but useful remark.) 
The proof is elementary, cf. [GMS2].

\subsection{Lemma}
 {\it Let $X$ be a compact K\"ahlerian manifold. For any $i\geq 0$ consider 
a part of the holomorphic de Rham complex
$$
\Omega^{[i}_X:\ 0 \lra \Omega^i_X \lra \Omega^{i+1}_X \lra \ldots
$$
with $\Omega^i_X$ in degree $0$. 

For each $n$ we have a canonical isomorphism
$$
H^n(X,\Omega^{[i}_X) \isom \oplus_{p = 0}^n H^{n-p}(X,\Omega_X^{i+p})
$$
(Hodge decomposition).}

This is the main assertion of the Hodge theory. 

\subsection{ }
 Using the previous two lemmas we get as a corollary a convenient description of a cohomology 
space playing an important role in the theory of chiral differential operators. 

Let $X$ be a smooth algebraic projective variety. We use an abbreviated notation from [GMS1]
$$
\Omega^{[2,3>}_X:\ 0\lra \Omega^2_X \lra \Omega_X^{3,cl} \lra 0
$$
Here the differential forms are algebraic. From the previous and GAGA we get a canonical isomorphism
$$
H^1(X,\Omega^{[2,3>}_X)\isom H^0(X,\Omega^3_X)\oplus H^1(X,\Omega^2_X) 
\eqno{(2.3.1)}
$$
Note that
$$
H^0(X,\Omega^{[2,3>}_X)\isom H^0(X,\Omega^2_X)
\eqno{(2.3.2)}
$$
for each smooth $X$.

\subsection{ }
From now on we return to the framework and notation of Introduction; in particular $X$ will be an abelian variety. 

We have canonical isomorphisms
$$
H^p(X,\Omega^q_X) \isom Hom(\Lambda^q\fg, \Lambda^p\hfg)
$$

\subsection{Some $\BC$-vertex algebroids }
 We shall use the terminology from [GMS1]. 

Suppose we are given an element $\lambda\in (\Lambda^3\fg)^* = Hom(\Lambda^3\fg,\BC)$. We can 
consider $\lambda$ as a skew symmetric function $\lambda(x_1,x_2,x_3),\ x_i\in \fg$. 

Define a $\BC$-vertex algebroid $A(\lambda)$ as follows. In the notation of {\it op. cit.}, 1.4, we set 
$$
A(\lambda) = (\BC, \fg, \fg^*, 0, 0, \langle , \rangle, c_\lambda)
$$
Here $\fg$ is considered as a Lie algebra with zero bracket, the action of $\fg$ on $\fg^*$ is trivial, 
the only nonzero components of the map 
$$
\langle , \rangle:\ (\fg\oplus\fg^*) \times (\fg\oplus\fg^*) \lra \BC
$$
are the obvious maps $\fg\times\fg^*\lra\BC$,  $\fg^*\times\fg\lra\BC$. The only nontrivial element of the definition is the map
$$
c_\lambda: \fg\times\fg \lra \fg^*
$$
defined by
$$
c_\lambda(x,y)(z) = \lambda(x,y,z)
$$

The only nontrivial axiom among (A1) - (A5) from {\it loc. cit.} is (A4) which is fulfilled due to 
complete skew symmetry of $\lambda$. 

Thus we get a set of vertex algebroids indexed by $(\Lambda^3\fg)^*$. Define a category 
$\CV alg(\fg)$ whose set of objects are these algebroids and morphisms are the morphisms of vertex 
algebroids in the sense of {\it op. cit.}, Thm. 3.5, of the form
$(\Id_\BC, \Id_\fg, \Id_{\fg^*}, h)$. 

Here $h: \fg \lra \fg^*$ must be skew-symmetric due to {\it loc. cit.} $(3.5\langle,\rangle)$. 

It follows that $\CV alg(\fg)$ is a groupoid with $Hom(A(\lambda),A(\lambda')) = \emptyset$ 
for $\lambda\neq \lambda'$ and 
$$
Hom(A(\lambda),A(\lambda)) = (\Lambda^2\fg)^*,
\eqno{(2.5.1)}
$$
the composition of morphisms being the addition in $(\Lambda^2\fg)^*$. 

\subsection{ }
Let $\CDO(X)$ denote the category of cdo's over $X$. The tangent bundle $\CT_X$ is trivial, 
so $ch_2(\CT_X) = 0$ and this category is nonempty. It is a groupoid whose set of isomorphism 
classes $\pi_0\CDO(X)$ is a torsor under $H^1(X,\Omega^{[2,3>}_X)$. 

In fact we can explicitely construct a big full subcategory of $\CDO(X)$. Identify 
$\fg$ with left invariant vector fields on $X$; it gives rise to a trivialisation 
$$
i:\ \fg\otimes\CO_X \iso \CT_X
$$
Using the push-forward construction [GMS1], 1.8, 1.9, we associate to each vertex algebroid $A(\lambda)$ 
from 2.5 an $\CO_X$-vertex algebroid, and hence a cdo, over $X$. Let us denote this cdo by 
$i_*A(\lambda) = \CA(\lambda)$ and call the cdo's of the form $\CA(\lambda)$ {\it quasi-constant}.  

In particular $\CA(0)$ is a distinguished object in $\CDO(X)$ which allows to identify 
$\pi_0\CDO(X)$ with 
$$
H^1(X,\Omega^{[2,3>}_X) \isom H^0(X,\Omega^3_X)\oplus H^1(X,\Omega^2_X) 
\isom (\Lambda^3\fg)^*\oplus Hom(\Lambda^2\fg,\hfg)
\eqno{(2.6.1)}
$$
The group of automorphisms of any object in $\CDO(X)$ is 
$$
H^0(X,\Omega^{[2,3>}_X) \isom H^0(X,\Omega^2_X) \isom (\Lambda^2\fg)^*
\eqno{(2.6.2)}
$$

Let 
$$
\CDO_c(X) = \CDO^{0,3}(X)\subset \CDO(X)
$$ 
denote the full subcategory spanned by quasi-constant cdo's. 
We have
$$
\pi_0\CDO^{0,3}(X)\isom H^0(X,\Omega^3_X)\isom (\Lambda^3\fg)^*
$$  
The push-forward defines an equivalence of groupoids
$$
i_*:\ \CV alg(\fg) \iso \CDO_c(X)
\eqno{(2.6.3)}
$$
Indeed, $i_*$ induces an isomorphism on $Hom$'s due to (2.6.2) and (2.5.1).

Thus, $\CDO(X)$ has an "easy to describe" part $\CDO^{0,3}(X)$ and "not easy to describe" part 
$\CDO^{1,2}(X)$ which is a full subcategory with
$$
\pi_0\CDO^{1,2}(X) \isom H^1(X,\Omega^2_X) 
\isom  Hom(\Lambda^2\fg,\hfg)
$$

\subsection{ }
Let $\CB$ be a non-degenerate tdo over $X$ with the characteristic class 
$$
\mu = \mu_\CB \in H^1(X,\Omega^1_X) = Hom(\fg,\hfg)
$$
We consider $\mu_\CB$ as an isomorphism $\fg\iso \hfg$. It induces isomorphisms
$$
(\Lambda^i\fg)^*\iso (\Lambda^i\hfg)^*,\ \lambda\mapsto\lambda_\mu
\eqno{(2.7.1)}
$$
sending $\lambda:\ \Lambda^i\fg\lra\BC$ to the composition 
$$
\Lambda^i\hfg\overset{\Lambda^i\mu^{-1}}\lra \Lambda^i\fg 
\overset{\lambda}\lra\BC
$$
This gives rise to an equivalence of categories
$$
F_\mu:\ \CV alg(\fg) \iso \CV alg(\hfg),\ F_\mu(A(\lambda)) = A(\lambda_\mu)
$$
More precisely, $F_\mu$ acts on objects via isomorphisms (2.7.1) with $i = 3$ and on morphisms 
via (2.7.1) with $i = 2$. 
 
It follows an equivalence
$$
F_\mu:\ \CDO_c(X) \iso \CDO_c(\hX),\ F_\mu(\CA(\lambda)) = \CA(\lambda_\mu), 
$$
cf. (2.6.3). 

This is our "Fourier - Mukai transform" on the easy to describe part of $\CDO(X)$.  

\subsection{ }
Let us extend $F_\mu$ to the whole of $\CDO(X)$. To do this we replace 
$\CDO(X)$ by an equivalent "easy to describe" groupoid $\CDO(X)^\sim$. 

Namely, the set of objects of $\CDO(X)^\sim$ will be by definition equal to 
$H^1(X,\Omega^{[2,3>}_X)$ and the group of automorphisms of each object equal to 
$H^0(X,\Omega^{[2,3>}_X)$. A functor
$$
\CDO(X) \lra \CDO(X)^\sim
\eqno{(2.8.1)}
$$
sending a cdo on $X$ to its isomorphism class is an equivalence of categories. 

Given a class $\mu$ as in 2.7, we define an equivalence of groupoids
$$
F_\mu:\ \CDO(X)^\sim \iso \CDO(\hX)^\sim
\eqno{(2.8.2)} 
$$
by the same token as above. Namely, taking into account the identifications (2.6.1) and (2.6.2), 
the action of $F_\mu$ on objects is induced by an isomorphism
$$
(\Lambda^3\fg)^*\oplus Hom(\Lambda^2\fg,\hfg) \iso 
(\Lambda^3\hfg)^*\oplus Hom(\Lambda^2\hfg,\fg)
$$
which is the direct some of an isomorphism 
$$
(\Lambda^3\fg)^* \iso (\Lambda^3\hfg)^*
$$
already defined in (2.7.1) and of an isomorphism
$$
Hom(\Lambda^2\fg,\hfg) \iso Hom(\Lambda^2\hfg,\fg)
$$
defined as follows. To an element $\nu\in Hom(\Lambda^2\fg,\hfg)$ we associate 
the composition 
$$
\Lambda^2\hfg \overset{\Lambda^2\mu^{-1}}\lra \Lambda^2\fg
\overset{\nu}\lra \hfg \overset{\mu^{-1}}\lra \fg
\eqno{(2.8.3)}
$$
The action of $F_\mu$ on the morphisms comes through an isomorphism
$$
(\Lambda^2\fg)^* \iso (\Lambda^2\hfg)^*
$$
defined in (2.7.1). 

This completes the definition of the equivalence (2.8.2). Using (2.8.1) we get finally an equivalence 
$$
F_\mu:\ \CDO(X) \iso \CDO(\hX)
\eqno{(2.8.4)} 
$$

\subsection{ }
 One can do the same thing with the usual tdo's. According to [BB], the groupoid 
$\TDO(X)$ of tdo's on $X$ has
$$
\pi_0\TDO(X) = H^1(X,\Omega^{[1,2>}_X) \isom 
$$
$$
\isom H^1(X,\Omega^1_X)\oplus H^0(X,\Omega^2_X),
$$
and for any $\CB\in Ob \TDO(X)$
$$
Aut(\CB) = H^0(X,\Omega^{[1,2>}_X) \isom H^0(X,\Omega^1_X)
$$
Using the same device as in 2.8 for each non-degenerate $\CB$ with $c(\CB) = \mu$ we get the equivalences
$$
F_\mu:\ \TDO(X)\iso \TDO(\hX)
$$
Although we have not verified this, it is natural to expect that 
$- F_\mu(\CB)$ is isomorphic to the Fourier - Mukai transform $\Phi(\CB)$ constructed in [PR], cf. (1.1.2).

\subsection{Relative version}
 Let $\pi: X \lra S$ be a smooth map. One can define a groupoid 
$\CDO(X/S)$ of cdo's over $S$; its objects will be  
vertex algebroids of the form $(\CO_X, \CT_{X/S}, \Omega_{X/S}, \ldots)$. 

If $X$ is an abelian veriety over $S$, we can repeat the previous construction to obtain functors 
$$
F_{S,\mu}:\ \CDO(X/S) \lra \CDO(\hat X/S)
$$
where
$$
\mu \in \Gamma(S,R^1\pi_*\Omega^1_{X/S})
$$
is a non-degenerate class, i.e. its restriction to each fiber $X_s, s\in S$ is non-degenerate. 

\subsection{Example: Hitchin fibration}
For all the terminology below see [DP] and references therein. Let $C$ be a smooth proper curve 
of genus $> 0$, $G$ a simple complex Lie group. Fix a maximal torus $T\subset G$; set $\ft = Lie(T)$.
 
Consider the Hitchin fibration
$$
h:\ \Higgs(G) \lra B(G)
$$
The total space $\Higgs(G)$ parametrizes $K_C$-valued $G$-Higgs bundles over $C$ where $K_C$ is the canonical 
class; it more or less coincides with the cotangent space to the space of (semistable) $G$-bundles on $C$.

The base 
$$
B(G) = H^0(C,(K_C\otimes\ft)/W)
$$
where $W$ is the Weyl group of $G$.   

We have $\pi_0(\Higgs(G)) = \pi_1(G)$; let $\Higgs_0(G)$ be the connected component corresponding to the 
identity in $\pi_1(G)$. 

Outside certain discriminantal divisor $\Delta(G) \subset B(G)$, i.e. over 
$$
B'(G) = B(G)\setminus \Delta(G),
$$ 
the fibers of $h$ are abelian varieties. So consider the restriction of $h$ 
$$
h'(G): \Higgs'_0(G) \lra B'(G)
$$

Let $^LG$ denote the Langlands dual group. In [DP] an isomorphism
$$
\ell_G:\ B(G) \iso B(^LG)
$$
is defined (the identity for simply laced $G$) which preserves the dicriminantal loci. The main result of {\it op. cit.} asserts that, after identification of the bases of Hitchin fibrations using $\ell_G$, 
$h'(G)$ and $h'(^LG)$ are dual abelian schemes. 

A scalar product $c$ on $\ft$ ("level") 
gives rise to a tdo $B_c$ on $\Higgs(G)$ which is fiberwise non-degenerate\footnote{we thank 
A.Beilinson for explaining this to us}. Thus we can apply the previous considerations to obtain the 
Fourier - Mukai transformation 
$$
F_c:\ \CDO(\Higgs'_0(G)/B'(G)) \iso \CDO(\Higgs'_0(^LG)/B'(^LG))
$$


\section{Quasiclassical considerations}



\subsection{The universal extension}
We keep the notations of the previous section. 
We have 
$$
Ext^1_{\CO_X}(\CO_X,\hfg^*_\CO) = H^1(X,\CO_X)\otimes\hfg^* = \hfg\otimes\hfg^*
$$
(here $?_\CO = ?\otimes\CO_X$), whence the canonical extension of 
$\CO_X$-modules 
$$
0 \lra \hfg^*_\CO \lra \CO_X^\natural \lra \CO_X \lra 0
\eqno{(3.1.1)}
$$
After passing to $\CO_X$-duals we get
$$
O \lra \CO_X\lra \CO_X^{\natural*} \lra\hfg_\CO \lra 0
$$ 

Consider a commutative $\CO_X$ algebra 
$$
\CA^\natural_X = Sym_{\CO_X}^\cdot(\CO_X^{\natural*})/(1_{\CO_X^\natural - 1_{Sym}})
$$
and let $X^\natural = \Spec \CA^\natural$. This is a commutative algebraic group, 
the universal extension of $X$ by a vector group:
$$
0 \lra \hfg^*  \lra X^\natural \lra X \lra 0  
$$

The FM transform of Rothstein, (cf. [R]), says that one has an equivalence of derived categories 
$$
D^b(\CD_X-mod) \iso D^b(\CO_{\hX^\natural}-mod),
\eqno{(3.1.2)} 
$$
so the algebra $\CD_X$ is FM dual to $\CA^\natural_{\hX}$. 

Now suppose we are given a tdo $\CD_\xi$ with the non-degenerate class
$$
\xi\in H^1(X,\Omega^1_X) = Hom(\fg^*,\hfg^*)
$$
It gives rise to 
$$
\xi^{-1}:\ \hfg^*_\CO \lra \fg^*_\CO = \Omega^1_X
\eqno{(3.1.3)}
$$
Taking the push forward of (3.1.1) by means of (3.1.2) we get an extension 
$$
0 \lra \Omega^1_X \lra E_\xi \lra \CO_X \lra 0
$$
whose class is $\xi$. 

So $E^*_\xi$ as an $\CO_X$-module is nothing but the the Picard algebroid corresponding to 
$\CD_\xi$:  
$$
E^*_\xi \isom \CD_{\xi,\leq 1} 
$$
\subsection{ }
On the other hand, consider the one-parametric family of tdo's 
$\CD_{t\xi}$, $t\in \BP^1$, as in [BB], 2.1.11. So "the quasi-classical limit" 
$\CD_{\infty\xi}$ is a Poisson algebra: the ring of functions on the twisted 
cotangent bundle $T^*_\xi$ on $X$. 

The underlying commutative ring is  
$$
Sym^\cdot \CD_{\xi,\leq 1}/(1_{\CD} - 1_{Sym}).
$$
It follows that $\xi$ induces an algebra isomorphism
$$
\xi_*:\ \CA_X^\natural \iso \CD_{\infty\xi}
\eqno{(3.2.1)}
$$
We get also a Poisson structure on $\CA_X^\natural$ coming from the one on 
$\CD_{\infty\xi}$.  

So we can say that 
$\CA_X^\natural$ is the limit of $\CD_{t\xi}$ when $t\ra\infty$ and 
$\CA_{\hX}^\natural$ is the limit of $\CD_{-t^{-1}\xi^{-1}}$ when $t\ra 0$. 

Here by $\CD_{-\xi^{-1}}$ we have denoted a tdo over $\hX$, the  FM transform of 
$\CD_\xi$.

Then the   
Rothstein equivalence (3.1.2) will be the limit of FM equivalences
$$
D^b(\CD_{t\xi}-mod) \iso D^b(\CD_{-t^{-1}\xi^{-1}}-mod)
$$
as $t\lra 0$.

\subsection{ }

Now suppose we are given a cdo $\CB$ over $X$. Applying to it the 
"FM in the direction $\xi$" we obtain a cdo $F_\xi(\CB)$ over 
$\hX$. 

The cdo's (as well as tdo's) form "a $\BC$-vector space in categories" and we have 
$$
F_{t\xi}(\CB) = t^{-1}F_{\xi}(\CB),\ t\neq 0.
$$
Passing to the limit $t\ra 0$ we get a quasiclassical object 
$$
(F_{\xi}(\CB))^{cl} = \lim_{t\ra 0} t^{-1}F_{\xi}(\CB)
$$
which is as a sheaf of  "twisted jets on $T^*\hX$" in the sense of 
[AKM].  

However, contrary to the case of tdo's, it seems to have no relation 
to the jets on $\hX^\natural$.

\bigskip\bigskip

\begin{center}
{\bf PART II: BOSON COMPACTFIED ON A TORUS}
\end{center}


\section{generalities}
\subsection{ }
\label{jetspaces}
Let $M,N$ be $C^\infty$-manifolds.  Denote by $J_\infty M_N$ the space of jets of sections of the trivial bundle $\pi: M_N\stackrel{\text{def}}{=}M\times N\rightarrow N$.
Reflecting the product structure of the jet space, the de Rham complex of $J_\infty M_N$ is bi-graded: $\Omega_{J_\infty M_N}^{\bullet}=\Omega_{M,N}^{\bullet\bullet}$.
In addition to the ordinary de Rham differential, it carries 2 supercommuting differentials: the vertical de Rham differential (one that is $\cO_N$-linear)
\[
\delta:\Omega_{M,N}^{\bullet\bullet}\rightarrow\Omega_{M,N}^{\bullet\bullet+1},
\]
and
\[
d:\Omega_{M,N}^{\bullet\bullet}\rightarrow\Omega_{M,N}^{\bullet+1\bullet}
\]
which is determined by the fact that $\cO_{J_\infty M_N}$ is a $D_N$-algebra.  In fact, $\Omega_{M,N}^{\bullet 0}$ is the de  Rham complex of $\cO_{J_\infty M_N}$ regarded as a $D_N$-module, and the extension to the entire $\cO_{J_\infty M_N}$ is determined by the condition $[d,\delta]=0$.

Consider the unit circle $S^1$ and a cylinder $\Sigma=S^1\times\RR$ equipped with coordinates, $\sigma$ and $(\sigma, \tau)$ resp., where $\tau$ is the
standard coordinate on $\RR$ and $\sigma$ is the angular coordinate on $S^1$.

The bicomplex we have just defined will arise as either $\Omega_{M,\Sigma}^{\bullet\bullet}$ or $\Omega_{T^*M,S^1}^{\bullet\bullet}$. Their various components
will carry various structure elements of the ``theory'' as follows:

\subsection{ }\label{lagr-and-gamma}
 A  Lagranagian $\cL$ is an element of $ \Omega_{M,\Sigma}^{20 }$.

 The differential of each Lagrangian can be uniquely written as
\begin{equation}
\label{var-1-form-arises}
\delta\cL=-d\gamma+\sum_i\frac{\delta\cL}{\delta x_i}\delta x_i,
\end{equation}
for some $\gamma\in \Omega_{M,\Sigma}^{11 }$ known as the {\em variational 1-form}.  Here we have  fixed a local coordinate system $\{x_1,x_2,...\}$ on $M$ although the term
$\sum_i\frac{\delta\cL}{\delta x_i}\delta x_i$ can be described without making any such choice.  In any case, the system $\delta\cL/\delta x_i=0$, $i=1,2,...$,
is known as the Euler-Lagrange equations or equations of motion. Let  $Sol\subset J_\infty M_\Sigma$ be defined by the equations  $\delta\cL/\delta x_i=0$, $i=1,2,...$, and all of their differential
consequences.

All of this is  a fancy way to write down the integration by parts  derivation of  the Euler-Lagrange equations  familiar to anybody from
a university calculus of variations course, but what that familiar derivation misses is the boundary terms, hence the variational 1-form $\gamma$, an important ingredient
of a field theory. 

\subsection{ }
\label{int-mot}
 {\bf Definition.} {\em An integral of motion is a 1-form $I\in\Omega_{M,\Sigma}^{10 }$ such that the restriction $dI|_{Sol}=0$.}

The meaning of this is clear:  A classical trajectory of the string is a map $\Phi:\Sigma\rightarrow M$. It naturally lifts to a map
$\Sigma\rightarrow J_\infty M_\Sigma$, to be also denoted by $\Phi$, so that $\Phi(\Sigma)\subset Sol$.  The position of the string at ``time'' $\tau$ is the restriction of this
map to $S^1\times\{\tau\}$. Any $\beta\in\Omega_{M,\Sigma}^{10 }$  defines a function
\[
\RR\rightarrow\RR,\; \tau\mapsto\int_{S^1\times\{\tau\}}\Phi^*\beta,
\]
If $\beta=I$ is an integral of motion, then this function is independent of $\tau$, as the following computation shows:
\[
\partial_\tau\int_{S^1\times\{\tau\}}\Phi^*I=\int_{S^1\times\{\tau\}}\Phi^*\partial_\tau I\stackrel{dI=0}{=}-\int_{S^1\times\{\tau\}}\Phi^*\partial_\sigma I=0.
\]

\subsection{ }
\label{noter} N\"other's theorem  attaches an integral of motion to a symmetry of $\cL$.

A vertical vector field $\hat{\xi}\in\cT_{J_\infty M_\Sigma}$ is called {\em evolutionary} if it respects the $\cD_\Sigma$-module structure.

As a practical matter, any infinitesimal diffemorphism
$x\mapsto x+\epsilon F$, $F\in\cO_{J_\infty M_\Sigma}$, uniquely extends to an evolutionary vector field by $\partial_{\vec{j} }x\mapsto\partial_{\vec{j} }  x+\epsilon \partial_{\vec{j} }F$.

{\bf Definition.} {\em An evolutionary vector field $\hat{\xi}$ is a symmetry of $\cL$ if $\hat{\xi}\cL=d\alpha_\xi$ for some $\alpha_\xi\in\Omega^{10}_{M,\Sigma}$.}

N\"other's theorem states that
\begin{equation}
\label{form-noter-thm}
\hat{\xi}\text{  is a symmetry of } \cL\Longrightarrow \alpha_{\xi}-\iota_{\hat{\xi}}\gamma\text{ is an integral of motion,}
\end{equation}
where $\gamma$ is the variational 1-form.
  Furthermore, the space of integrals of motion
carries a bracket so that the corresponding Lie algebra maps onto the Lie algebra of symmetries of $\cL$.

\subsection{ }
\label{coisson-on-jet}
Somewhat separately from the above, if $M$ is a Poisson manifold, then $\Omega^{10}_{M,S^1}$, or rather the push-forward onto $S^1$ ,$(\pi_\infty)_* \Omega^{10}_{M,S^1}$,
is a coisson algebra. The chiral bracket (see the introduction to the Beilinson-Drinfeld book) can be defined by the following  formula. Let $\{.,.\}$ be the Poisson bracket
on $\cO_M$.  Then set, somewhat symbolically, for $f,g\in\cO_M\subset \Omega^{10}_{M,S^1}$
\[
\{f(\sigma)d\sigma,g(\sigma')d\sigma'\}=\{f,g\}(\sigma')\delta(\sigma-\sigma').
\]
This uniquely extends to a bracket on the entire $\Omega^{10}_{M,S^1}$ for the reason that $\cO_{J_\infty M_\Sigma}$ is a universal $D_{S^1}$-algebra generated by
$\cO_{M_{S^1}}$.

This construction works just as well for $M$ replaced with any open subset $U\subset M$ and gives us a sheaf of coisson algebras over $M$; this is the reason for skipping
$(\pi_\infty)_* $ in the notation.

 If we replace $M$ with $T^*M$ equipped with the canonical
Poisson bracket, then locally there is a coordinate system $q_i$, $p_i$, $1\leq i\leq \text{dim}N$, so that $\{p_i,q_j\}=-\delta_{ij}$, and the chiral bracket on $\Omega^{10}_{M,S^1}$  is
determined by the assignment
\[
\{p_i(\sigma)d\sigma, q_j(\sigma')d\sigma'\}=-\delta_{ij}\delta(\sigma-\sigma'),
\]
from which formulas such as $\{ p_i(\sigma)d\sigma, \partial_{\sigma'}q_j(\sigma')d\sigma'\}=\delta_{ij}\partial_\sigma\delta(\sigma-\sigma')$ follow by definition. This gives, of course,
a quasiclassical limit of the $\beta\gamma$-system.

\subsection{ }
\label{fourier-comp}
Attached to any coisson algebra is a Lie algebra of ``Fourier components of fields.'' In our case, it is
$h(\Omega^{10}_{M,S^1})\stackrel{\text{def}}{=}\Omega^{10}_{M,S^1}/d\Omega^{00}_{M,S^1}$.  For example, if $f,g$ are functions on $S^1$, then the Beilinson-Drinfeld prescription gives the bracket
\begin{eqnarray}
\nonumber& &\{\int f(\sigma)p_i(\sigma)d\sigma, \int g(\sigma')\partial_{\sigma'}q_j(\sigma')d\sigma'\}_{\text{Lie}}=\int f(\sigma)g(\sigma')\delta_{ij}\partial_\sigma\delta(\sigma-\sigma')\\
\nonumber&=&\int f(\sigma')g(\sigma')\delta_{ij}\partial_\sigma\delta(\sigma-\sigma')+\delta_{ij}\int f'(\sigma')g(\sigma')(\sigma-\sigma')\partial_\sigma\delta(\sigma-\sigma')\\
\nonumber&=&-\delta_{ij}\int f'(\sigma)g(\sigma)d\sigma,
\end{eqnarray}
where $\int$ means the class of the indicated form modulo exact forms.  Note that the summand $\int f(\sigma')g(\sigma')\delta_{ij}\partial_\sigma\delta(\sigma-\sigma')$ vanishes because it lies in the image of $\partial_\sigma$.

If furthermore we denote $\alpha^{(i)}_m=\int e^{-im\sigma}p_i(\sigma)d\sigma$,
$\beta^{(j)}_n=\int e^{-in\sigma}\partial_\sigma q_j(\sigma)d\sigma$, then, $\int e^{i(m+n)}d\sigma$ being $\delta_{m,-n}d\sigma$, we obtain the Heisenberg algebra
relations
\[
\{\alpha^{(i)}_m,\beta^{(j)}_n\}_{\text{Lie}}=im\delta_{ij}\delta_{m,-n}d\sigma.
\]
\subsection{ }
\label{rel-cotang-to-int-mot}

The relation of the coisson algebra structure on $\Omega^{10}_{M,S^1}$ to the field theory determined by a Lagrangian $\cL\in\Omega^{20}_{M,\Sigma}$ as outlined in
sect.~\ref{lagr-and-gamma} -- \ref{noter} is as follows. It was shown in [F.Malikov ``Lagrangian approach...''] that for a class of Lagrangians there is a diffeomorphism
$Sol\stackrel{\sim}{\rightarrow} J_\infty T^*M_{S^1}\times\RR$. Introduce $Sol^o = J_\infty T^*M_{S^1}\times\{0\}\subset Sol$.
One can show ({\em loc. cit.})  that the composite restriction map
 \[
 \Omega^{10}_{M,\Sigma}\rightarrow \Omega^{10}_{Sol}\rightarrow\Omega^{10}_{Sol^o}\stackrel{\sim}{\rightarrow} \Omega^{10}_{T^*M,S^1}
 \]
  induces an embedding of the Lie algebra of integrals of motion mentioned in sect.~\ref{noter} in the Lie algebra $h( \Omega^{10}_{T^*M,S^1})$. Thus a single coisson algebra $\Omega^{10}_{T^*M,S^1}$ carries a variety of field theories and what distinguishes among them is a collection of integrals motion, especially the  Hamiltonian, realized as a subspace of $\Omega^{10}_{T^*M,S^1}$.  Let us see how all of this plays out in the case of the $\sigma$-model on a circle.

\section{Boson compactified on a circle}
\subsection{ }
\label{sett-up-model}
As the target space $M$ we choose a radius $R$ circle, which we realize as the quotient $\RR/L$, where $L$ is the lattice $\{2\pi Rm,m\in\ZZ\}$.  A typical $\sigma$-model Lagrangian
is
\begin{equation}
\label{real-lagr}
\cL^{re}=\frac{1}{2}(\partial_\tau x^2-\partial_\sigma x^2)d\tau\wedge d\sigma,
\end{equation}
where $x$ is the standard coordinate on $\RR$. It will be more convenient to replace $\tau$ with $-i\tau$ (the ``Wick rotation'') and focus on
\begin{equation}
\label{compl-lagr}
\cL=\frac{i}{2}(\partial_\tau x^2+\partial_\sigma x^2)d\tau\wedge d\sigma.
\end{equation}
Technically, this requires that we replace all the arising sheaves, such as $\Omega^{\bullet\bullet}_{M,\Sigma}$, with their complexifications, such as $\BC\otimes\Omega^{\bullet\bullet}_{M,\Sigma}$. We assume this replacement made but keep the notation unchanged.

We have
\begin{equation}
\label{var-of-lagr-circle}
\delta\cL=-d(i\partial_\tau x\delta x\wedge d\sigma-i\partial_\sigma x\delta x\wedge d\tau)-i(\partial_\tau^2x+\partial_\sigma^2x)\delta x\wedge d\tau\wedge d\sigma.
\end{equation}
Thus we obtain the variational 1-form
\begin{equation}
\label{var-1-form-cicrle}
\gamma=i\partial_\tau x\delta x\wedge d\sigma-i\partial_\sigma x\delta x\wedge d\tau,
\end{equation}
variational 2-form
\begin{equation}
\label{var-2-form-cicrle}
\omega=\delta\gamma=i\delta(\partial_\tau x)\wedge\delta x\wedge d\sigma-i\delta(\partial_\sigma x)\wedge\delta x\wedge d\tau,
\end{equation}
and the solution space
\begin{equation}
\label{sol-circle}
Sol=\{\partial_\tau^2x+\partial_\sigma^2x=0\}.
\end{equation}
We see that the ring of functions
\[\BC[Sol]=\BC[S^1\times S^1\times\RR][\partial^\epsilon_\tau\partial_\sigma^n x,\; n\geq 0,\epsilon=0,1,n+\epsilon\geq 1].
\]
The infinite dimensional manifold $Sol$ shares many formal properties with $J_\infty S^1_\Sigma$; in particular it carries an action $\cD_\Sigma$, for example, $\partial_t$
operates as follows
\[
\partial_\sigma^n x\mapsto \partial_\sigma^n x+\epsilon\partial_\tau\partial_\sigma^n x,\;
\partial_\tau\partial_\sigma^n x\mapsto \partial_\tau\partial_\sigma^n x-\epsilon\partial_\sigma^{n+2} x.
\]
This action is invariant w.r.t. parallel transport along $\RR$ which allows us to get rid of $\tau$ by introducing
\[
Sol^o=\{\tau=0\}\subset Sol.
\]
We have
\begin{equation}
\label{ring-of-fnctn}
\BC[Sol^o]=\BC[S^1\times S^1][\partial^\epsilon_\tau\partial_\sigma^n x,\; n\geq 0,\epsilon=0,1,n+\epsilon\geq 1].
\end{equation}
  This indeed allows us to identify
$Sol^o$ with $J_\infty T^*S^1_{S^1}$ with $x$ a coordinate on the target $S^1$ and $\partial_\tau x$ the ``dual'' coordinate. The situation at hand is very simple and this identification may
seem rather arbitrary, but notice that (\ref{var-2-form-cicrle}) is reminiscent of the standard symplectic form and suggests the coisson bracket discussed in sect.~\ref{coisson-on-jet}:
\begin{equation}
\label{coiss-br-bos-circle}
\{\partial_\tau x(\sigma)d\sigma,x(\sigma')d\sigma'\}=i\delta(\sigma-\sigma').
\end{equation}
The same can be done for an arbitrary target $M$ with Riemannian metric.

\subsection{ }
\label{int-mot-circl}
$\BC[J_\infty S^1_\Sigma]$ is a $\cD_\Sigma$-module. For example, the action of $\partial_\tau$ is given by
\[
\partial_\tau\mapsto\partial_\tau+\widehat{\partial_\tau},\;\widehat{\partial_\tau}=(\partial_\tau x)\frac{\delta}{\delta x}+(\partial_\tau \partial_\sigma x)\frac{\delta}{\delta \partial_\sigma x}+\cdots
\]
The evolutionary vector field $\widehat{\partial_\tau}$ is a symmetry of $\cL$, sect.~\ref{noter},  because $\cL$ does not explicitly depend on $\tau$. One has
\[
\widehat{\partial_\tau}\cL=d(\frac{i}{2}(\partial_\tau x^2+\partial_\sigma x^2) d\sigma),
\]
and so $\alpha_{\partial_\tau}=\frac{i}{2}(\partial_\tau x^2+\partial_\sigma x^2) d\sigma$. Furthermore,
\[
\iota_{\widehat{\partial_\tau}}\gamma=i(\partial_\tau x)^2d\sigma-i\partial_\tau x\partial\sigma xd\tau.
\]
According to (\ref{form-noter-thm}), the corresponding integral of motion, in fact the Hamiltonian, is
\begin{equation}
\label{hamilt-circle}
H_{\partial_\tau}=-\frac{i}{2}((\partial_\tau x)^2-(\partial_\sigma x)^2)d\sigma+i\partial_\tau x\partial_\sigma xd\tau.
\end{equation}
A similar computation for $\partial_\tau$ replaced with $\partial_\sigma$ gives another integral of motion
\begin{equation}
\label{d-sigma-circle}
H_{\partial_\sigma}=-i\partial_\tau x\partial_\sigma xd\sigma-\frac{i}{2}((\partial_\tau x)^2-(\partial_\sigma x)^2)d\tau.
\end{equation}

One also notes that the constant vector field $x\mapsto x+\epsilon$, to be denoted $\delta/\delta x$, preserves the Lagrangian, and the corresponding integral of motion
is the {\em momentum}:
\begin{equation}
\label{comp-of-momentu}
H_{\delta/\delta x}=-\iota_{\delta/\delta x}\gamma= -i\partial_\tau xd\sigma+i\partial_\sigma xd\tau.
\end{equation}

\subsection{ }
\label{how-hamilt-gen-dynamics}
We will now illustrate the discussion of sect.~\ref{rel-cotang-to-int-mot} and show an example of how restricting to $Sol^o$ facilitates the computation of brackets.

Restricting to $Sol^o$ essentially means setting $\tau=0$, $d\tau=0$; for example, $ H_{\partial_\tau}|_{Sol^o}=-\frac{i}{2}((\partial_\tau x)^2-(\partial_\sigma x)^2)d\sigma$.

Computing as in sect.~\ref{fourier-comp}, one verifies that indeed the evolutionary vector field $\widehat{\partial_\tau}$ preserves $Sol^o$, where it equals
$\{\int H_{\partial_\tau}|_{Sol^o},.\}_{\text{Lie}}$:
\[
\{\int H_{\partial_\tau}|_{Sol^o},.x\}_{\text{Lie}}=\{\int -\frac{i}{2}((\partial_\tau x)^2-(\partial_\sigma x)^2)d\sigma,x\}_{\text{Lie}}=\partial_\tau x,
\]
\[
\{\int H_{\partial_\tau}|_{Sol^o},.\partial_\tau x\}_{\text{Lie}}=\{\int -\frac{i}{2}((\partial_\tau x)^2-(\partial_\sigma x)^2)d\sigma,\partial_\tau x\}_{\text{Lie}}=\partial_\sigma x(\sigma)\delta '_\sigma(\sigma-\sigma')\equiv -\partial_\sigma^2x(\sigma)\text{mod}\partial_\sigma,
\]
as it should in view of the equations of motion (\ref{sol-circle}).
\subsection{ }
\label{conf-symmetry}
The computations we have just carried out are but the tip of an iceberg that is conformal symmetry. Introduce complex-valued functions $z=\tau+i\sigma$ and $\bar{z}$.
It is easy to check that
\[
\cL=-\partial_z x\partial_{\bar{z}}x dz\wedge d\bar{z}.
\]
It follows easily from this presentation that any holomorphic (or antiholomorphic) vector field $f(z)\partial_z$ ($g(\bar{z})\partial_{\bar{z}}$)
is a  symmetry of $\cL$.  For example,
\[
f(z)\partial_z\cL=-d(f(z)(\partial_z x )^2 d\bar{z}),
\]
and so,
\[
\alpha_{f(z)\partial_z}=-f(z)(\partial_z x )^2 d\bar{z}.
\]
The variational 1-form can be rewritten as
\[
\gamma=\partial_z x\delta x\wedge dz-\partial_{\bar{z}}x\delta x\wedge d\bar{z}.
\]
Computing $\alpha_{f(z)\partial_z}-\iota_{\widehat{f(z)\partial_z}}\gamma$, see (\ref{form-noter-thm}), yields the corresponding integrals of motion
\begin{equation}
\label{int-conf-symm}
H_{f(z)\partial_z}=-f(z)(\partial_z x)^2 dz.
\end{equation}
Similarly,
\begin{equation}
\label{int-conf-symm-antiholo}
H_{g(\bar{z})\partial{\bar{z}}}=g(\bar{z})(\partial_{\bar{z}} x)^2 d\bar{z}.
\end{equation}

\subsection{ }
\label{coiss-alg-str}
Building on sect.~\ref{how-hamilt-gen-dynamics}, we will now restrict  the various integrals of motion to $Sol^o$ and see what kind of relations we obtain inside the standard
$\Omega^{10}_{S^1, S^1}\stackrel{\sim}{\rightarrow} \Omega^{10}_{Sol^o}$.  We get
(committing  the abuse of keeping the notation unchanged), cf. sect.~\ref{int-mot-circl},
\[
H_{\partial_\tau}=-\frac{i}{2}((\partial_\tau x)^2-(\partial_\sigma x)^2)d\sigma,\; H_{f(z)\partial_z}=-if(z)(\partial_z x)^2 d\sigma,\;H_{g(\bar{z})\partial{\bar{z}}}=-ig(\bar{z})(\partial_{\bar{z}} x)^2 d\sigma.
\]
Note that
\begin{equation}
\label{hequals-l-plus-l}
H_{\partial_\tau}=H_{\partial_z}+H_{\partial_{\bar{z}}}.
\end{equation}
In addition, motivated by the integrals of motion $H_{f(z)\partial_z}$ and $H_{g(\bar{z})\partial{\bar{z}}}$, introduce
$i\partial_z x d\sigma$, $i\partial_{\bar{z}}x d\sigma$.

The indicated elements give the following chiral brackets, all immediate consequences of (\ref{coiss-br-bos-circle}),
\begin{equation}
\label{chir-brack-heis}
\{i\partial_z x(\sigma)d\sigma,i\partial_z x(\sigma')d\sigma'\}=\frac{1}{2}\partial_\sigma\delta(\sigma-\sigma'),\;
\{i\partial_{\bar{z}} x(\sigma)d\sigma,i\partial_{\bar{z}} x(\sigma')d\sigma'\}=-\frac{1}{2}\partial_\sigma\delta(\sigma-\sigma') ,
\end{equation}
\begin{equation}
\label{chir-br-vir+}
\{-i(\partial_z x)(\sigma)^2 d\sigma,-i(\partial_z x)(\sigma')^2 d\sigma'\}=2(\partial_z x)(\sigma')^2 \partial_\sigma\delta(\sigma-\sigma')-
\partial_{\sigma'}((\partial_z x)(\sigma')^2) \delta(\sigma-\sigma'),
\end{equation}
\begin{equation}
\label{chir-br-vir-}
\{-i(\partial_{\bar{z}} x)(\sigma)^2 d\sigma,-i(\partial_{\bar{z}} x)(\sigma')^2 d\sigma'\}=-2(\partial_{\bar{z}} x)(\sigma')^2 \partial_\sigma\delta(\sigma-\sigma')+ \partial_{\sigma'}
(i(\partial_{\bar{z}} x)(\sigma')^2) \delta(\sigma-\sigma').
\end{equation}
The computations similar to those in sect.~\ref{fourier-comp} will give us 2 commuting copies of the Heisenberg algebra, spanned by
$\int i e^{im\sigma}\partial_z x(\sigma)d\sigma$ and  $\int ie^{im\sigma}\partial_{\bar{z}} x(\sigma)d\sigma$ (resp.), and (centerless) Virasoro algebra, spanned
by $\int(-1) i e^{im\sigma}(\partial_z x)(\sigma)^2 d\sigma$ and  $\int(-1) ie^{im\sigma}i(\partial_{\bar{z}} x)(\sigma)^2 d\sigma$ (resp.):
\begin{equation}
\label{lie-brack-heis+}
\{\int ie^{im\sigma}\partial_z x(\sigma)d\sigma,\int ie^{in\sigma}\partial_z x(\sigma)d\sigma\}=-\frac{1}{2}m\delta_{m,-n}id\sigma,
\end{equation}
\begin{equation}
\label{lie-brack-heis-}
\{\int ie^{im\sigma}\partial_{\bar{z}} x(\sigma)d\sigma,\int ie^{in\sigma}\partial_{\bar{z}} x(\sigma)d\sigma\}=\frac{1}{2}m\delta_{m,-n}id\sigma,
\end{equation}
\begin{equation}
\label{lie-br-vir+}
\{\int(-1)ie^{im\sigma}(\partial_z x)(\sigma)^2 d\sigma,\int(-1)ie^{in\sigma}(\partial_z x)(\sigma')^2 d\sigma\}=(m-n)\int(-1)ie^{i(m+n)\sigma}(\partial_z x)(\sigma)^2 d\sigma,
\end{equation}
\begin{equation}
\label{lie-br-vir-}
\{\int(-1)ie^{im\sigma}(\partial_{\bar{z}} x)(\sigma)^2 d\sigma,\int(-1)ie^{in\sigma}(\partial_{\bar{z}} x)(\sigma')^2 d\sigma\}=-(m-n)\int(-1)ie^{i(m+n)\sigma}(\partial_{\bar{z}x})(\sigma)^2 d\sigma.
\end{equation}

 Denote by $\cH^+$, $\cH^-$ and $\cV ir^+$, $\cV ir^-$ the respective Lie algebras.  Note another group of familiar brackets
 \begin{equation}
 \label{a+-vir+}
 \{\int(-1)(\partial_z x)(\sigma)^2 d\sigma,\int(-1)ie^{in\sigma}\partial_z x)(\sigma) d\sigma\}=-\frac{1}{2}n\int(-1)ie^{in\sigma}\partial_z x)(\sigma) d\sigma.
 \end{equation}
\begin{equation}
 \label{a--vir-}
 \{\int(-1)(\partial_{\bar{z}} x)(\sigma)^2 d\sigma,\int(-1)ie^{in\sigma}\partial_{\bar{z}} x)(\sigma) d\sigma\}=\frac{1}{2}n\int(-1)ie^{in\sigma}\partial_{\bar{z}} x)(\sigma) d\sigma.
 \end{equation}
The signs are different, this is quite crucial.

\subsection{ }
\label{quantization-beginning}
So far all our considerations have been quasiclassical. How does all of this quantize?  It is of course true that  $\Omega^{10}_{T^*S^1, S^1}$, where most of our computations have taken place, is a quasiclassical limit of a CDO $\cD_{S^1}$. Taking the latter for a quantization, however, is utterly unphysical. 

In most general terms, a quantum theory is a pair, an algebra and an appropriate module over this algebra referred to as the space of states.  Whatever the word ``algebra''
means in the present context, the most one can hope for is to have (a completion of) the universal enveloping algebra of the Lie algebra $\cH$, the Lie subalgebra of
$h(\Omega^{10}_{T^*S^1, S^1})$ generated by 
\[
\int e^{im\sigma}\partial_\tau x(\sigma) d\sigma\text{ and }\int e^{im\sigma}x(\sigma) d\sigma;
\]
 this is because the entire $h(\Omega^{10}_{T^*S^1, S^1})$ is generated
in a sense by the indicated elements. In other words, we wish to understand which elements of $\cH$ quantize. The answer depends on what we want from the ``module.''
There are three requirements.

(1) It is natural to demand that each of its elements (each state, a physicist would say) is a linear combination  of fixed energy states.  In other words,
the Hamiltonian, $\int H_{\partial_\tau }$, must be diagonalizable.

(2) A physicist requires  that the operator $\int H_{\partial_\tau}$ have a locally bounded from below spectrum, i.e., that any element generate a submodule
with spectrum of $H_{\partial_\tau}$ bounded from below.

(3) The spectra of the elements $\int \partial_\sigma x d\sigma$ and $i\int \partial_\tau x d\sigma$ satisfy
\begin{equation}
\label{spec-x-sigma}
\text{Spec}(\int \partial_\sigma x d\sigma)=L,
\end{equation}
\begin{equation}
\label{spec-x-tau}
\text{Spec}(i\int \partial_\tau x d\sigma)=L^*,
\end{equation}
where $L=\{2\pi Rm,\;m\in\ZZ\}$, the lattice via which our $S^1$ was realized, see sect.~\ref{sett-up-model}, and $L^*=\{(2\pi R)^{-1}m,\;m\in\ZZ\}$ the dual lattice. This is the first
instance of geometry of the target entering the fray, the preceding discussion being independent of the radius of the circle and equally applicable to the target $\RR$.

The meaning of condition (\ref{spec-x-sigma}) is rather clear: $x$ being the coordinate on the covering space, $\int \partial_\sigma x d\sigma$ attaches to a string position
$x: S^1\rightarrow S^1$ its length. The latter can be thought of as $2\pi R\times$the number of times the string wraps around the origin, known as  the {\em   winding  number}.

As to (\ref{spec-x-tau}), observe that $i\int \partial_\tau x d\sigma$ is interpreted as the average {\em momentum} , see (\ref{comp-of-momentu}), and so (\ref{spec-x-tau}) becomes a familiar quantum mechanical assertion that for a particle on a circle
the momentum quantizes on the dual lattice.

Therefore the joint spectrum is very simple:
\[
\text{Spec}(\int \partial_\sigma x d\sigma,i\int \partial_\tau x d\sigma)=L\oplus L^*.
\]

\bigskip

It is then immediate to obtain
\begin{equation}
\label{spec-in-terms-of left-right-heis}
\text{Spec}(\int i \partial_z x d\sigma,-\int i\partial_{\bar{z}} x d\sigma)=\{\frac{1}{2}(l-l^*,l+l^*),\; l\in L, l^*\in L^*\}
\end{equation}

\subsection{ }
\label{phys-math-transl-alg}

Condition (1) at least suggests (if not determines)  the ``algebra.''  Indeed, $\cH$ is closed under the action of $\{\int H_{\partial_\tau },.\}$ and  we denote by
$\cH^{ss}\subset\cH$ the linear span of all eigenvectors of $\{\int H_{\partial_\tau },.\}$. It follows that in any reasonable module  each eigenvector of $H_{\partial_\tau}$ generates an $\cH^{ss}$-submodule with a diagonalizable action of $H_{\partial_\tau}$. It is then natural to expect that in order to have (1) satisfied the module should carry an action
of $\cH^{ss}$, not the entire $\cH$.

Note the commutation relations, which follow from (\ref{hequals-l-plus-l}) and (\ref{a+-vir+}-\ref{a--vir-}),
\begin{equation}
\label{eigenv-a+}
\{\int H_{\partial_\tau },\int ie^{im\sigma}\partial_z x(\sigma)d\sigma\}=-\frac{1}{2}m\int ie^{im\sigma}\partial_z x(\sigma)d\sigma,
\end{equation}
\begin{equation}
\label{eigenv-a-}
\{\int H_{\partial_\tau },\int ie^{im\sigma}\partial_{\bar{z}} x(\sigma)d\sigma\}=\frac{1}{2}m\int ie^{im\sigma}\partial_{\bar{z}} x(\sigma)d\sigma,
\end{equation}
It follows that $\cH^+\oplus\cH^-\subset\cH^{ss}$, and it is easy to understand that, in fact, $\cH^+\oplus\cH^-=\cH^{ss}$. (Note that the only basis element left out is
$\int x(\sigma)d\sigma$, and it generates a $2\times 2$ Jordan cell.)

The Heisenberg $\cH^+\oplus\cH^-$ is then tentatively assigned to be the algebra. Let us now find the module.

\subsection{ }
\label{phys-math-transl-modul}

A typical $\cH^+ $-module is the Fock module $V^+_{a}=\text{Ind}_{\cH^+_+}^{\cH^+}\BC_a$, where $\cH^+_+\subset\cH^+$ is the subalgebra spanned by
$\int e^{im\sigma}\partial_z xd\sigma$, $m\geq 0$, and  $\BC_a$ is its 1-dimensional representation, where $\int e^{im\sigma}\partial_z xd\sigma$, $m> 0$, act as 0 and the central element
$\int \partial_zd\sigma$ as multiplication by $a$.

Similarly defined  is the $\cH^- $-module  $V^-_{a}=\text{Ind}_{\cH^-_+}^{\cH^-}\BC_a$, where $\cH^-_+\subset\cH^+$ is the subalgebra spanned by
$\int e^{im\sigma}\partial_{\bar{z}} xd\sigma$, $m\geq 0$.

One can as well induce from the opposite subalgebra, such as $\cH^+_-\subset\cH^+$, the subalgebra spanned by
$\int e^{im\sigma}\partial_z xd\sigma$, $m\leq 0$.  The result can be equivalently described as the restricted dual. Thus we obtain modules $(V^+_{a})^*$ and
$(V^-_{a})^*$.

It follows from (\ref{eigenv-a+},\ref{eigenv-a-}) that the requirements sect.~\ref{quantization-beginning}(1,2) combined
imply that the space of states has a filtration by $\cH^+\oplus\cH^-$-modules $V^+_{a^+}\otimes (V^-_{a^-})^*$, $(a^+,a^-)\in\BC$; this is where the different signs in
(\ref{a+-vir+},\ref{a--vir-}) show up.

In fact, one can prove that under very mild restrictions the actual space of states is a direct sum of the indicated
representations. 

Finally, the requirement sect.~\ref{quantization-beginning}(3) determines which highest weights $(a^+,a^-)$ do occur. As (\ref{spec-in-terms-of left-right-heis}) implies,
he answer is
\begin{equation}
\label{determ-state-space}
\text{the space of states }=\bigoplus_{l\in L,l^*\in L^*}V^+_{1/2(l-l^*)}\otimes (V^-_{-1/2(l+l^*)})^*.
\end{equation}
\subsection{ }
\label{vert-alg-str}
The indicated space of states carries a rich multiplicative structure. First of all $V^+_0$ is a vertex algebra.  In order to set up this structure we need to work over a disc rather than
a cylinder. We embed the latter into $\BC$ with coordinate $v$ by letting $v=e^z$.  Denote $\alpha_m=\int ie^{im\sigma}\partial_z x(\sigma) d\sigma$ and attach the formal operator-valued
series $\alpha(v)=\sum_m v^{-m-1}\alpha_m$ to $\alpha_{-1}1\in V^+_0$. It is well known that this assignment uniquely extends to a vertex algebra structure on $V^+_0$.

It follows from (\ref{lie-brack-heis+}) that
\[
\alpha(u)\alpha(v)\sim\frac{ -1/2}{(u-v)^2}.
\]
Furthermore, the action of the algebra of holomorphic vector fields, see (\ref{int-conf-symm}), is generated by the Virasoro field $-:\alpha(u)^2:$.

A similar discussion applies to $(V^-_0)^*$, but in this case we face a dilemma.  Since we need to deal with the dual module, which is turned upside down, so to say, we replace
$e^{i\sigma}$ with $e^{-i\sigma}$, and so we have to decide whether we use $v^{-1}|_{S^1}=e^{-i\sigma}$ or $\bar{v}|_{S^1}=e^{-i\sigma}$.  There are at least two reasons why
$\bar{v}$ is the right choice. First, the reason to introduce $v=e^z$ is to compactify the cylinder by adding the origin $v=0$, where $v^{-1}$ is not defined. Second,
it is  the algebra of antiholomorphic vector fields (\ref{int-conf-symm-antiholo}) that operates on $(V^-_0)^*$.

Therefore,  we define $\bar{\alpha}_m=\int ie^{-im\sigma}\partial_{\bar{z}}x(\sigma) d\sigma$,
$\bar{\alpha}(\bar{u})=\sum_m \bar{u}^{-m-1}\alpha_m$. We have,
\[
\bar{\alpha}(\bar{u})\bar{\alpha}(\bar{v})\sim\frac{ -1/2}{(\bar{u}-\bar{v})^2},
\]
and  (\ref{int-conf-symm-antiholo}) becomes the Virasoro field $-:\alpha(\bar{u})^2:$.

All these definitions are made so as to have the Lie brackets of sect.~\ref{coiss-alg-str} preserved.

\subsection{ }
\label{vertex-operators}
$V^+_0\otimes (V^-_0)^*$ is then the ``algebra,'' but  the entire space of states also carries a multiplicative structure of sorts. This is about  a classic and very familiar
construction attaching a vertex operator to elements of Heisenberg algebra modules. In order to recall it, it is convenient to introduce some notation.

Define an inner product
on $\BC\alpha$ and $\BC\bar{\alpha}$ by setting $(\alpha,\alpha)=-1/2=(\bar{\alpha},\bar{\alpha})$. This identifies these spaces with their duals and we can write, for example,
$V^+_{m\alpha}$ instead of $V^+_{-1/2m}$.

Next, consider the group algebra $\BC[\BC\alpha]$ and denote its typical element by $e^{m\alpha}$, $m\in\BC$. The canonical generator
$1_{m\alpha}=1\in\BC=\BC_{m\alpha}\subset
\text{Ind}_{\cH^+_+}^{\cH^+}\BC_{m\alpha}$ will be identified with $e^{m\alpha}$.

Attached to each $a\in V^+_{m\alpha}$ is
a formal series $Y(a,u)$ s.t. for each $k$ and  $b\in V^+_{k\alpha}$ one has
\[
Y(a,u)b\in u^{-1/2mk}V^+_{(k+m)\alpha}((u)).
\]
For example, $Y(e^{m\alpha},u)$ is the familiar {\em vertex operator} $e^{m\alpha}(u)$ defined by
\[
e^{m\alpha}(u)=e^{m\alpha}u^{m\alpha_0}:\exp{(-\sum_{n=-\infty}^{+\infty}\frac{nm\alpha}{n}u^{-n})}:.
\]
It is  straightforward (and classic) to compute (see the Kac book)
\[
e^{m\alpha}(u)e^{n\alpha}(v)=(u-v)^{-1/2mn}:e^{m\alpha}(u)e^{n\alpha}(v):.
\]
Here is what this means.  Fix $m_j$, $b_j\in V^+_{m_j\alpha}$, $0\leq j\leq n$. Then the formal series
\begin{eqnarray}
\nonumber
& &Y(b_n,u_n)Y(b_{n-1},u_{n-1})\cdots Y(b_1,u_1)b_0\\
&\in& V^+_{\sum_jm_j\alpha}[[u_1,...,u_n]][...u_j^{-1},...,(u_i-u_j)^{-1}...]\prod_{0<i<j}(u_i-u_j)^{-1/2m_im_j}\prod_{i>0}u_i^{-1/2m_im_0},\nonumber
\end{eqnarray}
where appropriate Taylor  series expansions are made.
It follows (cf.  the Frenkel-Ben Zvi book) that the (obviously defined) correlators
\[
\langle Y(b_n,u_n)Y(b_{n-1},u_{n-1})\cdots Y(b_1,u_1)\rangle
\]
are multi-valued functions  on $(\BC^!)^n$ without diagonals behaving as $(u_i-u_j)^{-1/2m_im_j}$ near the diagonal $\{u_i=u_j\}$.

Were all $m_j$ even numbers, we would thus obtain a vertex algebra structure. 

Of course a similar discussion applies to the antiholomorphic part of the story, and we obtain vertex operators, such as $e^{m\bar{\alpha}}(\bar{u})$ or, more generally,
$Y(\bar{b},\bar{u})$, $b\in(V^-_{-m\bar{\alpha}})^*$.

Now let us combine the two in the case of the space of states (\ref{determ-state-space}). Attached to each $b\otimes \bar{b}\in V^+_{1/2(l-l^*)}\otimes (V^-_{-1/2(l+l^*)})^*$
is a vertex operator $Y(b\otimes \bar{b};u,\bar{u})=Y(b,u)\otimes Y(\bar{b},\bar{u})$.  Throwing in another $a\otimes \bar{a}\in V^+_{1/2(m-m^*)}\otimes (V^-_{-1/2(m+m^*)})^*$,
one obtains the expansion
\begin{eqnarray}
& &Y(b\otimes \bar{b};u,\bar{u})Y(a\otimes \bar{a};v,\bar{v})\nonumber\\
&=&(u-v)^{-1/2(l-l^*)(m-m^*)}(\overline{u-v})^{-1/2(l+l^*)(m+m^*)}:Y(b\otimes \bar{b};u,\bar{u})Y(a\otimes \bar{a};v,\bar{v}):.\nonumber
\end{eqnarray}
Now one is pleased to notice that thanks to the nontrivial choice of the holomorphic/antiholomorphic highest weights an {\em a priori} multi-valued function acquires a canonical
univalued branch. This kind of multiplicative structure was axiomatized by Kapustin and Orlov  [KO] and called... vertex algebra, having reserved the name of chiral algebra for what
mathematicians tend to call a vertex algebra.

\subsection{T-duaiity.}
\label{t-duaity}
Consider the involution on the set of radii: $R\mapsto (4\pi^2R)^{-1}$. The corresponding models we constructed are obviously isomorphic -- simply because
under this involution $L$ and $L^*$ are swapped. This can also be expressed as an unexpected duality between momentum modes and winding modes, see
sect.~\ref{quantization-beginning}.  The T-duality is a prototype of the mirror symmetry.

But what if a model and its dual are the same, i.e., $R=(4\pi^2R)^{-1}$?

\subsection{Chiral algebra.}
\label{chiral-a;gebra}
Define the chiral algebra of model  (\ref{determ-state-space}) to be the ``centralizer'' of the antiholomorphic Virasoro algebra, i.e., $\cap_{n\geq -1}\text{Ker}(-\bar{\alpha}(\bar{u})^2)_{n}$.

It is easy to see that the answer is
\[
\bigoplus_{l=-l^*}V^+_{2l\alpha}\otimes e^0.
\]
Therefore, generically the answer is the uninspiring vertex algebra $V^+_0$. At the other extreme is the case where $L^*=L$, which occurs precisely at the fixed point
$R=(4\pi^2R)^{-1}$. Then we obtain the chiral algebra  equal to\footnote{we are indebted to A.Kapustin who explained this to us, [Kap1]}
\[
\bigoplus_{l}V^+_{2l\alpha},
\]
which is the lattice vertex algebra on an even lattice: $\cV_{2\alpha\ZZ}$.

 \subsection{ }
 \label{disc-alg-vs-statespace}
 One often hears that the chiral algebra defines a theory.   This is  wrong for the torus model we have just discussed: if $L$ is generic, then the chiral algebra is just a Heisenberg
 algebra, it does not know about the lattice.
 
 One often hears, too, that a lattice vertex algebra on an even lattice is a $\sigma$-model on the corresponding torus. This is  true only with various qualifications.
 In our vase we can make a precise statement.

 If $L=L^*$, then the arising $\cV_{2\alpha\ZZ}$ has precisely 2 modules, itself and another one, which  looks like an ``algebra'' on the ``lattice'' $(2\ZZ+1)\alpha$.
 So denote it by $\cV_{(2\ZZ+1)\alpha}$. Since $l-l^*\equiv l+l^*\text{ mod }2$, it appears that (\ref{determ-state-space}) can be rewritten \`a la WZW:
 \[
\cV_{2\alpha\ZZ}\otimes \bar{\cV}_{2\alpha\ZZ}\oplus\cV_{(2\ZZ+1)\alpha}\otimes \bar{\cV}_{(2\ZZ+1)\alpha}.
\]

\section {higher dimensional torus}
\label{higher dimensional torus}
The higher dimensional situation is very similar, the only conceptually new ingredient being a ``$B$-field,'' and we will be brief.

\subsection{ } 
\label{set-up-multidim}
We consider $\RR^n$ with a positive symmetric bilinear form $g=g(.,.)$ and a rank $n$ lattice $L\subset \RR^n$. Define $M=\RR^n/L$.

In addition to $g(.,.)$, which will also be regarded as a metric on the tangent bundle to $M$, we fix an antisymmetric bilinear form on $\RR^n$, $B=B(.,.)$, to be likewise regarded
as a differential 2-form on $M$.

Since we have the canonical  coordinate system $x=(x^1,...,x^n)$, we obtain matrices $g_{ij}=(\partial_i,\partial_j)$ and $b_{ij}=B(\partial_i,\partial_j)$.

Expressions such as $g(\partial_\tau x,\partial_\tau x)$ or $B(\partial_\tau x,\partial_\sigma x)$ can be  thought of as the value of the indicated form on the
indicated ``generic'' tangent vectors, but in reality they are simply the quadratic functions on the jet space, $g(\partial_\tau x,\partial_\tau x)=g_{ij}\partial_\tau x^i\partial_\tau x^j$
and
$B(\partial_\tau x,\partial_\sigma x)=b_{ij}\partial_\tau x^i\partial_\sigma x^j$.

We will also regard $g$ and $B$ as linear transformations $\RR^n\rightarrow (\RR^n)^*$ defined by $\langle g(x),y\rangle=g(x,y)$ and
$\langle B(x),y\rangle=B(x,y)$. Expressions, such as $g(\partial_\tau x)$ and $B(\partial_\sigma x)$ will be interpreted similarly.

\subsection{ }
\label{lagr-and-such-multidim}
We consider
\[
\cL=\frac{\sqrt{-1}}{2}(g(\partial_\tau x,\partial_\tau x)+g(\partial_\sigma x,\partial_\sigma x)-2\sqrt{-1}B(\partial_\tau x,\partial_\sigma x))d\tau\wedge d\sigma.
\]
The equations of motion and the variational 1-form are
\[
\partial_\tau ^2 x +\partial_\sigma^2 x=0,
\]
\[
\gamma=(\sqrt{-1}g(\partial_\tau x,\delta x)+B(\partial_\sigma x,\delta x))d\sigma-(\sqrt{-1}g(\partial_\sigma x,\delta x)+B(\partial_\tau x,\delta x))d\tau.
\]
The former equation shows that $Sol^o$ is canonically $J_\infty TM_{S^1}$, the latter
 allows us to identify $Sol^o$ with $J_\infty T^*M_{S^1}$ via a version of the Legendre transform, where the canonical coordinates on $T^*M$ are $\{x^j\}$ and are chosen to be
\begin{equation}
\label{canon-coord-jigh-dim}
p_j^B=\sqrt{-1}g(\partial_\tau x)_j+ B(\partial_\sigma x)_j\stackrel{\text{def}}{=}\sqrt{-1}g_{ij}\partial_\tau x^i+b_{ij}\partial_\sigma x^i.
\end{equation}
The ``canonical commutation relations'' are read off the formula for $\gamma$:
\begin{equation}
\label{caonn-comm-rel-high-dim}
\{p^B_i(\sigma),x^j(\sigma')\}=\delta_i^j\delta(\sigma-\sigma'),
\end{equation}
and so the identification $Sol^o\stackrel{\sim}{\rightarrow}J_\infty T^*M_{S^1}$ is that of manifolds with coisson structure.

Note that the apparent dependence of the coisson structure on $B$ is inessential:  for any closed 2-form $\alpha$, the map 
\[
p_j\mapsto p_j+\iota_{\partial_j}\alpha\stackrel{\text{def}}{=}p_j+\alpha_{ji}\partial_\sigma x^i
\]
defines an automorphism of the coisson algebra $J_\infty T^*M_{S^1}$.

\subsection{ }
\label{int-mot-high-dim}
Corresponding to the symmetries of $\cL$, $\partial_j=\delta/\delta x^j$, $\partial_\tau$, $f(z)\partial/\partial z$, $g(\bar{z})\partial/\partial \bar{z}$, are the integrals of motion:

momentum
\begin{equation}
\label{mom-high-dim}
H_{\delta/\delta x^j}=-p^B_j,
\end{equation}
Hamiltonian
\begin{equation}
\label{Hamil-high-dim}
H_{\partial_\tau}=-\frac{\sqrt{-1}}{2}(g(\partial_\tau x,\partial_\tau x)-g(\partial_\sigma x,\partial_\sigma x))d\sigma,
\end{equation}
and 2 commuting copies of the Virasoro algebra:
\begin{equation}
\label{VIRA-high-dim}
H_{f(z)\partial/\partial z}=-if(z)g(\partial_z x,\partial_z x)d\sigma,
\end{equation}
\begin{equation}
\label{antiVIRA-high-dim}
H_{f(\bar{z})\partial/\partial \bar{z}}=-if(\bar{z})g(\partial_{\bar{z}} x,\partial_{\bar{z}} x)d\sigma,
\end{equation}
Of course
\[
H_{\partial_\tau}=H_{\partial/\partial z}+H_{\partial/\partial \bar{z}}.
\]
Of these formulas, only the first one involves $B(.,.)$.
\subsection{ }
\label{hamilt-eigenv-quantiz-high-dim}
As we have seen, the prerequisite for quantization is the diagonalization of the Hamiltonian $\int H_{\partial_\tau}$ as an operator acting on
$h(\Omega^{10}_{T^*M_{S^1}})$.  It follows from sect.~\ref{int-mot-high-dim} that the answer is given by the Lie algebras $\cH^+$ and $\cH^{-}$
spanned by $\{\int e^{im\sigma}\partial_z x^j d\sigma\}$ and $\{\int e^{im\sigma}\partial_{\bar{z}} x^j d\sigma\}$ resp. The Lie brackets follow from the
Poisson brackets
\begin{eqnarray}
\{\sqrt{-1}\partial_z x^i(\sigma)d\sigma,\sqrt{-1}\partial_z x^j(\sigma')d\sigma'\}&=&\frac{1}{2}g^{ij}\partial_\sigma\delta(\sigma-\sigma'),\nonumber\\
\{\sqrt{-1}\partial_{\bar{z}} x^i(\sigma)d\sigma,\sqrt{-1}\partial_{\bar{z}} x^j(\sigma')d\sigma'\}&=&-\frac{1}{2}g^{ij}\partial_\sigma\delta(\sigma-\sigma'),\nonumber\\
\{\partial_zx^i(\sigma),\partial_{\bar{z}}x^j(\sigma')\}&=&0,\nonumber
\end{eqnarray}
all consequences of (\ref{caonn-comm-rel-high-dim}), and imply that $\cH^+$ and $\cH^-$ are 2 commuting copies of the Heisenberg algebra ; here $(g^{ij})$ is the bilinear form on the cotangent bundle, the transfer of the original $(g_{ij})$ defined on the tangent bundle.

Exactly as in sect.~\ref{phys-math-transl-modul}, we define 
an $\cH^+ $-module  $V^+_{a}=\text{Ind}_{\cH^+_+}^{\cH^+}\BC_a$ and an $\cH^- $-module  $V^-_{a}=\text{Ind}_{\cH^-_+}^{\cH^-}\BC_a$, where $a$ is an element of
$\BC^n$. The latter space appears here as the  complexification of the fiber of $TM=T\RR^n/L$ and has the meaning of the dual to the center of $\cH^+$  spanned
by $\{\int\partial_z x^j d\sigma, j=1,...,n\}$ or $\{\int\partial_{\bar{z}} x^j d\sigma, j=1,...,n\}$.

Arguing as in sect.~\ref{int-mot-high-dim} we obtain
\begin{equation}
\label{space-states-high-dim-tentat}
\text{the space of states }=\bigoplus_{(a,b)\subset\Lambda}V^+_{a}\otimes (V^-_{b})^*.
\end{equation}
It remains to find $\Lambda\subset\BC^n$, the joint spectrum of $\{\int\partial_z x^j d\sigma, \int\partial_{\bar{z}} x^j d\sigma, j=1,...,n\}$.

\subsection{ }
\label{comp-spectrum-high-dim}
Arguing as in sect.~\ref{quantization-beginning}, we obtain a 1-1 correspondence between the joint spectrum of elements  $\{\int\partial_\sigma x^i d\sigma,\int p_j d\sigma\}$ and
$L\oplus L^*$ so that  $(l,l^*)\in L\oplus L^*$ defines an assignment
\[\int\partial_\sigma x^i d\sigma\mapsto \langle l,dx^i\rangle,\; \int p^B_j d\sigma\mapsto \langle l^*,\partial_j\rangle.
\]
More relevant for our purpose is $\int \partial_\tau x^j d\sigma$, dual to $\int p^B_j d\sigma$, and we similarly obtain, at the same spectrum point,
\[
\sqrt{-1}\int \partial_\tau x^j d\sigma\mapsto \langle g^{-1}(l^*)-g^{-1}\circ B(l), dx^j\rangle,
\]
where $g^{-1}: (\RR^n)^*\rightarrow \RR^n$ and $B: \RR^n\rightarrow (\RR^n)^*$ are the maps naturally induced by the metric $(.,.)$ and the 2-form $B$ resp; this follows from
(\ref{canon-coord-jigh-dim}).

It follows that the points of $\Lambda$ are parametrized by $L\oplus L^*$ so that $(l,l^*)$
determines the functional
\begin{equation}
\label{joint-spectr-2heis-high-dim}
(\sqrt{-1}\partial_z x^j d\sigma,\sqrt{-1}\partial_{\bar{z}} x^j d\sigma)\mapsto (\frac{1}{2}\langle g^{-1}(l^*-B(l))-l, dx^j\rangle,\frac{1}{2}\langle g^{-1}(l^*-B(l))+l, dx^j\rangle)
\end{equation}

This along with (\ref{space-states-high-dim-tentat}) gives us the desired quantization. We shall now formulate the result using the standard lattice vertex algebra conventions.

\subsection{ }
\label{vert-alg-high-dim}
As discussed in sect.~\ref{lagr-and-such-multidim}, $Sol^o=J_\infty TM_{S^1}$ and contains 2 copies of the fiber of $T^*M$, spanned by $\partial_\tau x^j$ and
$\partial_\sigma x^j$ resp.  Upon complexification, we obtain $T^*_{p,\BC}M\oplus T^*_{p,\BC}M\subset\BC\otimes_{\RR} \cO_{J_\infty TM_{S^1}}$ and an operator of complex
conjugation, $\bar{ }$ .  We then have 2 more copies of $T^*_{p,\BC}M$, one spanned by $\partial_z x^j$, another by $\partial_{\bar{z}} x^j$, one being another's complex
conjugate. Since we have a canonical coordinate system $\{x^j\}$, we shall denote the former by $(\BC^n)^*$ and the later by $(\bar{\BC}^n)^*$.

Both  $(\BC^n)^*$ and $(\bar{\BC}^n)^*$ carry the bilinear form $-1/2g^{-1}=-1/2(g^{ij})$. We obtain 2 copies of the Heisenberg chiral algebra defined by the OPEs:
\[
\alpha(u)\beta(v)\sim-\frac{g^{-1}(\alpha,\beta)}{2(u-v)^2}, \;
\bar{\alpha}(\bar{u})\bar{\beta}(\bar{v})\sim-\frac{g^{-1}(\bar{\alpha},\bar{\beta})}{2(\bar{u}-\bar{v})^2},
\; \alpha,\beta\in(\BC^n)^*
\]
The Heisenberg vertex algebra modules will still be denoted by $V^+_a$ (or $V^-_a$), but  now it will be assumed that $a\in(\BC^n)^*$ and the value of, say,
$\alpha_{0}$ on the generator $e^a$ will be determined by the bilinear form:
\[
\alpha_{0}e^a=-1/2g^{-1}(\alpha,a)e^a.
\]

With all these conventions in place we summarize by stating that
\begin{equation}
\label{space-states-high-dim}
\text{the space of states }=
\bigoplus_{(l,l^*)\subset L\oplus L^*}V^+_{-l^*+B(l)+g(l)}\otimes V^-_{-l^*+B(l)-g(l)}.
\end{equation}
 The multiplicative structure on this space (state-field correspondence) is defined in the usual manner, reviewed briefly in sect.~\ref{vertex-operators}, the result being
 a Kapustin-Orlov type of vertex algebra.
 
 \subsection{ }
 At least when $B=0$, exactly as in the 1-dimensional case, one obtains
 
 $\bullet$ $T$-duality, i.e., an isomorphism  of thus defined vertex algebras  attached to the pairs $(g(L),L^*)$ and $(L^*,g(L))$ resp.;
 
 \bigskip
 
 $\bullet$ the chiral algebra of the model, which in the case where $L^*=g(L)$ is isomorphic to the lattice vertex algebra attached to the lattice $2L$.

\section{the baby wess-zumino-witten model  and cdo}
\label{the baby wess-zumino-witten model }
Our discussion can be generalized in at least 2 different ways\footnote{This subject goes somewhat beyond our original intentions and we will be brief}. 

\subsection{ }
\label{non-const-metr}
First, neither the metric $g(.,.)$ nor the 2-form $B(.,.)$ have to be constant; in fact,  $B(.,.)$ does not
have to be closed either. If so, it is neither natural nor necessary to assume that $M$ is a torus.

 Given any Lagrangian of the form
\[
\cL=\frac{\sqrt{-1}}{2}(g(\partial_\tau x,\partial_\tau x)+g(\partial_\sigma x,\partial_\sigma x)-2\sqrt{-1}B(\partial_\tau x,\partial_\sigma x))d\tau\wedge d\sigma,
\]
one still obtains the conformal symmetry manifesting itself in 2 Virasoro elements, 
\begin{equation}
\label{virasoro-general-nonconst}
-ig(\partial_z x,\partial_z x)d\sigma\text{  and }-ig(\partial_{\bar{z}} x,\partial_{\bar{z}} x)d\sigma,
\end{equation}
and the Hamiltonian equal to the sum
 \begin{equation}
 \label{hamilt-as-sum-wzw}
H_{\partial_\tau} -i\int (g(\partial_z x,\partial_z x)+g(\partial_{\bar{z}} x,\partial_{\bar{z}} x))d\sigma,
 \end{equation}
   just as in (\ref{VIRA-high-dim}, \ref{antiVIRA-high-dim}).
Note that none of this involves $B(.,.)$. Diagonalizing the Hamiltonian, however,  becomes problematic, as there is no obvious analogue of commuting pairs of bosons,
$\partial_z x^j$ and $\partial_{\bar{z}} x^j$.

\subsection{ }
\label{h-flux}
Another possibility is to introduce what physicists call an $H$-flux.  Fix a closed 3-form $H$ on $M$. Pick a fine enough open cover $\{U_j\}$ and a collection of
of 2-forms $B_j$ on $U_j$  so that  $dB_j=H$ on $U_j$ and $B_i-B_j$ is exact on $U_i\cap U_j$.  Given these data, we define $\cL^{H}$ to be a collection of
Lagrangians $\{\cL_j\}$, each $\cL_j$ defined over (jets in) $U_j$ by a familiar formula
\[
\cL_j=\frac{\sqrt{-1}}{2}(g(\partial_\tau x,\partial_\tau x)+g(\partial_\sigma x,\partial_\sigma x)-2\sqrt{-1}B_j(\partial_\tau x,\partial_\sigma x))d\tau\wedge d\sigma.
\]
Thanks to fact that $\cL_i-\cL_j$ is exact, much of the above carries over to this case; see [Mal] for details.  For example,  canonical commutation
relations
\[
\{p^{B_k}_i(\sigma),x^j(\sigma')\}=\delta_i^j\delta(\sigma-\sigma'),
\]
cf. (\ref{caonn-comm-rel-high-dim}), valid only over $U_k$ (where  $p^{B_k}_i$ is defined), are replaced with the globally defined
\begin{equation}
\label{twisted-sym-cdo}
\{p_i(\sigma),x^j(\sigma')\}=\delta_i^j\delta(\sigma-\sigma'),\;\{p_i(\sigma),p_j(\sigma')\}=\iota_{\partial_i}\iota_{\partial_j}H(\sigma')\delta(\sigma-\sigma'),
\end{equation}
where $p_i=\sqrt{-1}g(\partial_\tau x)_i$ and $\iota_{\partial_i}\iota_{\partial_j}H=h_{ijk}\partial_\sigma x^k$ provided
$H=h_{ijk}dx^i\wedge dx^j\wedge dx^k$.

Therefore  $\cL^H$ defines a coisson algebra structure on $J_\infty T^*M_{S^1}$ different from the canonical one that we have used so far. This fits with
the classification of CDOs  [GMS1]  recalled in sect.~2.6.  Namely, an obvious morphism of complexes 
\[
\Omega^{3,cl}_M\rightarrow(0\rightarrow\Omega^2_M\rightarrow\Omega^{3,cl}_M\rightarrow 0)
\]
induces a map
\[
H^0(M,\Omega^{3,cl}_M)\rightarrow H^1(M,\Omega^2_M\rightarrow\Omega^{3,cl}_M).
\]
Since we work with quasiclassical objects (rather than quantum ones, as the case was in Part I), equivalence classes of coisson algebras
are actually identified with the latter cohomology group (and are not only a torsor over the latter.)  The algebra that corresponds to the image of $H\in H^0(M,\Omega^{3,cl}_M)$ under this map is the one that we obtained via (\ref{twisted-sym-cdo}).

Note that exactly because the Virasoro elements \ref{virasoro-general-nonconst} are independent of $B$, a model with $H$-flux is conformally symmetric at the quasiclassical level.
Diagonalization of the Hamiltonian, however, is again problematic, and so is quantization. In fact, it is known that a (quantum) $\sigma$-model with
constant metric and  $H$-flux on a torus is not conformally
symmetric [Kap2].

\subsection{WZW}
\label{twzw}
It is all the more remarkable then that there is an important case where the 2 problems we have just discussed ( caused by a nonconstant metric or $H$-flux ) happen to be each
other's cure.

Let $G$ be a compact simple Lie group.  Define by $(.,.)$ a (unique up to proportionality)
bilinear form on $\fg= T_eG$. Then $([.,.],.)$ is   an invariant trilinear form on $\fg$. Let $g(.,.)$ and $H$ denote the corresponding
$G$-invariant metric and 3-form on $G$ For any $k\in\BC$, consider $\cL^{kH}$, the standard $\sigma$-model Lagrangian attached to $g(.,.)$ and $H$ in sect.~\ref{h-flux}.

It is clear  that  the action of $\fg$ on $G$ on the left and on the right defines symmetries of $\cL^{k/2H}$. Denote by 
\[
j_l^{k}, j_r^{k}:\fg\rightarrow \cO_{J_\infty T^*G_{S^1}},
\]
the corresponding morphisms. The morphisms $j_l^0$ and $ j_r^0$ are nothing but a composite of two canonical maps
\[
\fg\stackrel{\text{action}}\longrightarrow\cT_{G}\hookrightarrow \cO_{J_\infty T^*G_{S^1}},
\]
but they seem to be unrelated to the Virasoro algebras (\ref{virasoro-general-nonconst}) and the Hamiltonian of the model. 

The deformed maps are
\begin{equation}
\label{deformd-affine-alg}
j_l^k(x)=j_l^0(x)+\frac{k}{2}g(x,.),\;j_r^k(x)=j_l^0(x)-\frac{k}{2}g(x,.).
\end{equation}
These give rise to two coisson algebra embeddings
\begin{equation}
\label{chiral-deformd-alg-emb-aff}
j^k_l: V(\fg)^{coiss}_{k(.,.)}\rightarrow\cO_{J_\infty T^*G_{S^1}}\leftarrow V(\fg)^{coiss}_{-k(.,.)}: j^k_r,
\end{equation}
whose images centralize each other; here $V(\fg)^{coiss}_{k(.,.)}$ is a quasiclassical limit of the vertex algebra attached to $\fg$ with central charge $k(.,.)$.

Now notice that $V(\fg)^{coiss}_{k(.,.)}$ has a Virasoro element of its own, $L$. One can then verify that precisely when $k=2$, the following holds
\begin{equation}
\label{coinc-2-vir-deformed}
j^2_l(L)=-ig(\partial_z x,\partial_z x)d\sigma,\; j^2_r(L)=ig(\partial_{\bar{z}} x,\partial_{\bar{z}} x)d\sigma.
\end{equation}
Notice that the sign difference is easy to understand: while $j^2_l(L)$ and $j^2_r(L)$ are lifts of $i\partial_\sigma$, $-ig(\partial_z x,\partial_z x)d\sigma$ and
$-ig(\partial_{\bar{z}} x,\partial_{\bar{z}} x)d\sigma$ are lifts of $\partial_z$ and $\partial_{\bar{z}}$ resp.; since $\partial_z=1/2(\partial_\tau-i\partial_\sigma)$
and $\partial_{\bar{z}}=1/2(\partial_\tau+i\partial_\sigma)$, the former equals negative the latter when restricted to functions of $\sigma$ only.

It follows, cf. (\ref{hamilt-as-sum-wzw}),
\begin{equation}
\label{why-differ-sign-hamilt-baby-wzw}
H_{\partial_\tau}=\int j^2_l(L)-\int j^2_l(L)
\end{equation}

Therefore, thus arising 2 copies of the affine Lie algebra, one at level $2(.,.)$ another at level $-2(.,.)$,provide the desired diagonalization of the Hamiltonian. The computation
of the spectrum of momentum operators, which constituted an important  of the torus model analysis, has no place here, since $G$ is not an affine space.
Instead, one uses the $\fg\oplus\fg$-symmetry of the model to  demand that the space of states be a direct sum of finite dimensional $\fg\oplus\fg$-modules. Exactly how these modules are to be arranged is apparently a matter
of choice. A natural possibility is to use  $\cO_{J_\infty T^*G_{S^1}}$ itself as a guide. One has a  $V(\fg)^{coiss}_{k(.,.)}\otimes V(\fg)^{coiss}_{-k(.,.)}$-module decomposition
\begin{equation}
\label{chiral-peter-weyl}
\cO_{J_\infty T^*G_{S^1}}=\bigoplus V^{coiss}_{\lambda,2(.,.)}\otimes V^{coiss}_{\lambda^*,-2(.,.)},
\end{equation} 
where $V_{\lambda,k(.,.)}$ is the Weyl module induced from the simple finite dimensional $\fg$-module $V_\lambda$ at level $k(.,.)$ and $V^{coiss}_{\lambda,k(.,.)}$
is its quasiclassical limit.

Now we notice, as we have done twice previously, sect.~\ref{phys-math-transl-modul} and \ref{hamilt-eigenv-quantiz-high-dim} , that the boundedness of the energy
 from below and the minus sign in (\ref{why-differ-sign-hamilt-baby-wzw}) require that the right action  part be ``turned upside down, which leads
to a change of the sign of the level and gives
\begin{equation}
\label{def-wzw-quant-s[ace-state}
\text{WZW space of states}=\bigoplus V_{\lambda,2(.,.)}\otimes V_{\lambda^*,2(.,.)}
\end{equation}
at  level $(.,.)\neq0$. The same argument suggests that the chiral algebras generated by $j$

\subsection{Remarks}
\label{Remark}
The embeddings (\ref{chiral-deformd-alg-emb-aff}) are a rather simple quasiclassical limit of an appealing  analog valid at the level of chiral differential operator
algebras:
\[
j^k_l: V(\fg)_{k(.,.)}\rightarrow\cD^{ch}_{G,k}\leftarrow V(\fg)_{-(k+\kappa)(.,.)}: j^k_r,
\]
where $\kappa$ is a constant depending on the choice of $(.,.)$; if $(.,.)$ is the Killing form, then $\kappa$ is 1. This fact was observed in [FP]. A chiral algebra
proof can be found in [AG]; see also [GMS2]. It was related to the WZW model in  the quasiclassical limit in [Mal]. The reader will  notice that an attempt to quantize WZW
as an ``algebra'' leads to the quantum shift $-k\mapsto-k-\kappa$ and a loss of the antiholomorphic part.  
Physics quantization changes the sign, $-k\mapsto k$, retains the antiholomorphic part, but largely looses the algebra structure, which is replaced by the KZ equations
correlation functions, etc.

The important decomposition (\ref{chiral-peter-weyl}), or rather its more difficult quantum version, was proved in [FS].

\bigskip\bigskip

\centerline{\bf Bibliography}

\bigskip\bigskip

[AG] Arkhipov, S.; Gaitsgory, D. Differential operators on the loop group via chiral algebras, Int. Math. Res. Not. 2002, no. 4, 165Ð210

[AKM] T.Arakawa, T.Kuwabara, F.Malikov, Localisation of affine $W$-algebras, arXiv:1112.0089. 

[BB] A.Beilinson, J.Bernstein, A proof of Jantzen conjectures, in: I.M.Gelfand Seminar, Adv. in Sov. Math., 
{\bf 16}, Part 1 (1993), 1 - 50. 

[BD] A.Beilinson, V.Drinfeld, Chiral algebras, Colloquium Publications, vol.51, 2004

[DP] R.Donagi, T.Pantev, Langlands duality for Hitchin systems, arXiv:math/060461.  

[FP] Feigin, B. Parkhomenko, S. Regular representation of affine Kac-Moody alge- bras; in: Algebraic and geometric methods in mathematical physics (Kaciveli, 1993), 415 - 424, Math. Phys. Stud., 19, Kluwer Acad. Publ., Dordrecht, 1996.

[FS] Frenkel, I., Styrkas K., Modified regular representations of affine and Virasoro algebras, VOA structure and semi-infinite cohomology, posted on the archive: QA/0409117

[GMS1] V.Gorbounov, F.Malikov, V.Schechtman, Gerbes of chiral differential operators. II, {\it Inv. Math.} 
{\bf 155} (2004), 605 - 680. 

[GMS2] V.Gorbounov, F.Malikov, V.Schechtman, On chiral differential operators over homogeneous spaces, 
{\it Int. J. Math. Sci.} {\bf 26} (2001), 83 - 106. 

[Kap1] A.Kapustin, private communication.

[Kap2] A.Kapustin, private communication.

[KO]  A.Kapustin, D.Orlov, Vertex Algebras, Mirror Symmetry, And D-Branes: The Case Of Complex Tori, Comm. Math. Phys. 233 (2003), no. 1, 79--136 

[LL] D.L.Landau, E.M.Lifshitz, Quantum mechanichs, Course of Theoretical Physics, v.3, 3ed., Elsevier Science Ltd.

[Mal]  F.Malikov, Lagrangian approach to sheaves of vertex algebras,
{\it Comm.Math.Phys.}, {\bf 278} (2008), p.487-548

[M] D.Mumford, Abelian varieties, Bombay 1968. 

[P] A.Polishchuk, Abelian varieties, theta functions and the Fourier transform, Cambridge University Press. 

[PR] A.Polishchuk, M.Rothstein, Fourier transform for $D$-algebras, 
\newline arXiv:math/9901009. 

[R] M.Rothstein, Sheaves with connection on abelian varieties, alg-geom/9602023.

\bigskip

F.M.: Department of Mathematics, University of Southern California, Los Angeles, CA 90089, USA

V.Sch.: Institut de Math\'ematiques de Toulouse, Universit\'e Paul Sabatier, 31062 Toulouse, France

\end{document}